\documentclass[a4paper,10pt,openany]{article}
\usepackage[french,greek,english]{babel}
\usepackage[latin1]{inputenc}
\usepackage{color,listings,titlesec}
\usepackage[headings]{fullpage}
\usepackage{amsfonts}
\usepackage{fancyhdr}
\usepackage{graphicx}
\usepackage{amsmath}
\newcommand{\g}{\mathfrak{g}}

\newcommand{\h}{\mathfrak{h}}

\begin{document}

\thispagestyle{empty}
\begin{center}
\vspace{0.5cm}

{\large \textbf{Equivalence theorems and an explicit formula for the biquantization character.}}
\vspace{0.5cm}

Panagiotis Batakidis\footnote{Department of Mathematics, Aristotle University of Thessaloniki. e-mails: panagiotis.batakidis@gmail.com, batakidis@math.jussieu.fr}
\end{center}

\textbf{Abstract.} Let $\g$ be a Lie algebra, $\h\subset\g$ a subalgebra and $\lambda$ a character of $\h$ and $\epsilon$ a formal deformation parameter. In this paper we compute a character for the reduction algebra $H^0_{(\epsilon)}(\h_{\lambda}^{\bot},d^{(\epsilon)}_{\h_{\lambda}^{\bot},\mathfrak{q}})$ and its specialization $H^0_{(\epsilon=1)}(\h_{\lambda}^{\bot},d^{(\epsilon=1)}_{\h_{\lambda}^{\bot},\mathfrak{q}})$ via deformation quantization as in  \cite{K},\cite{CT},\cite{CF3}. Then we compare this construction with the computation of the character for the Penney eigendistribution $<\alpha(f),\psi>:=\int_{H/H\cap B}\psi(h)e^{f(H)}\mathrm{d}_{H/H\cap B}(h),\;\psi\in C^{\infty}_c(G,B,\chi_f)$, as in \cite{CG1},\cite{CG3},\cite{FUJI},\cite{FUJI2}. Using a double induction on $\dim\g,\;\dim\h$ it is proved that if the lagrangian condition ($\g$ nilpotent and $\exists\mathcal{O}\subset  \lambda+\h^{\bot},\;\textlatin{s.t}\;\forall l\in\mathcal{O},\; 2[dim (H~\cdot~l)] =dim (G\cdot l)$) holds, then the two computations coincide on the level of $H^0_{(\epsilon=1)}(\h_{\lambda}^{\bot},d^{(\epsilon=1)}_{\h_{\lambda}^{\bot},\mathfrak{q}})$. Finally we apply our results in an example of a 5-dimensional nilpotent Lie algebra and compute a term or order 3 with rational coefficients in the character formula.


\section{Introduction.}
In this section we recall  basic results from non-commutative harmonic analysis and review the construction of a character of $(U(\g)/U(\g)\h_f)^{\h}$ from representation theory. We also recall the facts and techniques that we need from deformation quantization. Our main references here are ~\cite{B1},~\cite{CG2},~\cite{FUJI},~\cite{FUJI2} and ~\cite{FUJI3} for the representation theoretic part, and \cite{K},\cite{CF2},\cite{CF3},\cite{CT},\cite{BAT1} for deformation quantization.

\textbf{General setting.} Let $G$ be a real nilpotent, connected and simply connected Lie group with $\g$ its Lie algebra, $\h$ a subalgebra of $\g$, $\lambda\in\h^{\ast}$ such that $\lambda([\h,\h])=0$, a character of $\h$. For $\;Y\in\h$, $\chi_{\lambda}:\;H\longrightarrow \mathbb{C}$ defined by $\chi_{\lambda}(\exp Y)=e^{i\lambda(Y)}$, is a unitary character of $H$, the Lie group associated to $\h$. Let $C^{\infty}(G,H,\chi_\lambda)$ be the vector space of complex smooth functions $\theta$ on $G$ that satisfy the property $\forall h\in H, \forall g\in G,\;\;\theta (gh)=\chi^{-1}_{\lambda}(h)\theta(g)$. Let also $\mathbb{D}(\g,\h,\lambda)$ be the algebra of linear differential operators, that leave the space $C^{\infty}(G,H,\chi_\lambda)$ invariant and commute with the left translation on $G$, i.e $\forall g\in G, \;\forall D\in \mathbb{D}(\g,\h,\lambda), \;\;\forall\theta \in C^{\infty}(G,H,\chi_\lambda)$, it is $D(C^{\infty}(G,H,\chi_\lambda))\subset C^{\infty}(G,H,\chi_\lambda),\;\;\textlatin{and}\;\;D(L(g)\theta)=L(g)(D(\theta))$. Koornwinder in \cite{koorn} proved that there is an algebra isomorphism
\begin{equation}\label{Koornwider}
(U(\g)/U(\g)\h_{\lambda})^{\h}\stackrel{\sim}{\rightarrow} \mathbb{D}(\g,\h,\lambda).
\end{equation} 
One can furthermore define the induced representation from $H,\chi_{\lambda}$ as $\tau_{\lambda}:=Ind(G\uparrow H,\chi_{\lambda})$ with Hilbert space $\mathcal{H}_{\lambda}:=L^2(G,H,\lambda)$ the separable completion of $C^{\infty}_c(G,H,\chi_\lambda)$, the compactly supported functions of $C^{\infty}(G,H,\chi_\lambda)$, with respect to the norm $||\phi||_2=\int_{G/H}|\phi(g)|^2d_{G/H}(g)$. The action of $G$ on $\phi\in L^2(G,H,\lambda)$ is translations by left: $\tau_{\lambda}(g)(\phi)(g')=\phi(g^{-1}g')$. These data correspond to a line bundle $\mathcal{L}_{\lambda}$ with base space $G/H$ and space of sections these functions $\phi$. Set $\g_{\mathbb{C}}:=\g\otimes\mathbb{C}$ and $U_{\mathbb{C}}(\g):=U(\g)\otimes\mathbb{C}$. Let $u\in(U_{\mathbb{C}}(\g)/U_{\mathbb{C}}(\g)\h_{i\lambda})^{\h}$, where $U_{\mathbb{C}}(\g)\h_{i\lambda}$ is the left ideal of $U_{\mathbb{C}}(\g)$ generated by $\h_{i\lambda}:=\{H + i\lambda(H),\;H\in\h\}$, and $D_u$ be the left invariant differental operator associated to $u$. The action of $\g$ on $\mathcal{H}_{\lambda}$ will be noted by $\mathrm{d}\tau_{\lambda}$ and for $X\in\g,\;g\in G$, it is defined by $\mathrm{d}\tau_{\lambda}(X)(\phi)(g)=\frac{d}{dt}\phi(\exp -tX\cdot g )|_{t=0}$. Sometimes we'll write also $\mathrm{d}\tau_{\lambda}(X)(\phi)(g)=L_X(\phi)(g)$. This action induces an action of $U_{\mathbb{C}}(\g)$, on $\mathcal{H}_{\lambda}$ that will be denoted also by $\mathrm{d}\tau_{\lambda}$. Let $\mathcal{H}^{-\infty}_{\lambda}$ be the space of antilinear continuous forms on $\mathcal{H}_{\lambda}^{\infty}$, the later being the dense subspace space of $C^{\infty}-$ vectors of $\mathcal{H}_{\lambda}$ (a function $\phi\in\mathcal{H}_{\lambda}$ is $C^{\infty}$ if the map $G\ni g\mapsto \tau_{\lambda}(g)\phi$ is smooth w.r.t the norm). The action of $U_{\mathbb{C}}(\g)$ on $\mathcal{H}^{-\infty}_{\lambda}$ will be denoted by $\mathrm{d}\tau_{\lambda}^{-\infty}$ and the action of $U_{\mathbb{C}}(\g)$ on $\mathcal{H}_{\lambda}^{\infty}$ is denoted respectively by $\mathrm{d}\tau_{\lambda}^{\infty}$. If $u=X_1X_2\cdots X_n,\;\;X_i\in\g$ is an element of $U_{\mathbb{C}}(\g)$ we set $u^t:=(-X_n)(-X_{n-1})\cdots(-X_1)$, $u^{\ast}:=(-\overline{X}_n)(-\overline{X}_{n-1})\cdots(-\overline{X}_1)$ where the bar stands for the complex conjugate. Let now $\{X_1,\ldots,X_n\}$ be a PBW basis of $\g$. For $u\in U_{\mathbb{C}}(\g)$, $\exists (k_s)_{s\in\mathbb{N}^n}$, $k_s\in\mathbb{C}$, such that $u$ can be expressed uniquely as $u=\sum_{s=(s_1,\ldots,s_n)} k_sX_1^{s_1}\cdots X_n^{s_n}$. Let $R(X)(\phi)(g):=\frac{d}{dt}{\phi(g\cdot\exp tX)}|_{t=0}$. Then the differential operator $D_u$ acts on the space of sections by $D_u(\phi)(g):= R_u(\phi)(g)$ where $R_u=\sum_s k_sR(X_1)^{s_1}\cdots R(X_n)^{s_n}$. Let $\mathfrak{b}$ be a polarization of $f\in\lambda+\h^{\bot}$ and $B$ the associated  to $\mathfrak{b}$ Lie subgroup of $G$. The map $\chi_f:\;B\longrightarrow \mathbb{C}$, defined as $\forall K\in\mathfrak{b},\;\;\chi_f(\exp K)=e^{if(K)},$ is a unitary character for $B$ and the induced representation from $B,\chi_f$ is $\tau_f=Ind(G\uparrow B,\chi_f)$ with Hilbert space $\mathcal{H}_f=L^2(G,B,f)$. 
The norm used here is the one coming from the product $\forall \phi_1,\forall \phi_2\in L^2(G,B,\chi_f),\;<\phi_1(x),\phi_2(x)>:=\int_{G/B}\phi_1(g)\overline{\phi_2(g)}\mathrm{d}_{G/B}(g)$, for an invariant measure $\mathrm{d}_{G/B}$ on $G/B$. Finally denote as $\mathcal{H}_f^{-\infty}$ the antidual space of the $C^{\infty}-$ vectors of $\mathcal{H}_f$ and as $\mathcal{H}_f^{-\infty,H}$ the $H-$ semi-invariant distribution vectors (\cite{FUJI2}).

\textbf{Penney vectors.} Let $\mathrm{d}_{H/H\cap B}$ be a left-invariant measure on $H/H\cap B$ and $\alpha_f\in\mathcal{H}_f^{-\infty}$ defined as

\begin{equation}\label{distribution}
<\alpha_f,\phi>=\int_{H/H\cap B}\overline{\phi(h)\chi_{\lambda}(h)}\mathrm{d}_{H/H\cap B}(h)\;\;\;\textlatin{for}\;\;\; \phi\in\mathcal{H}_f^{\infty}.
\end{equation}

The vector $\alpha_f$ is $H-$ semi-invariant (~\cite{B1}). Because of this equivariance property of $\alpha_f$, the algebra $\left(U_{\mathbb{C}}(\g)/U_{\mathbb{C}}(\g)\h_{-if}\right)^{\h}$ is acting on $\alpha_f$. Denote as $\h^{\bot}:=\{l\in\g^{\ast}/ l(\h)=0\}$ the annihilator of $\h$. A very important condition for this paper is the \textsl{lagrangian condition} which in algebraic terms is the following:
\begin{equation}\label{lagnow}
\exists \mathcal{O}\subset \lambda+\h^{\bot}\;\;\textlatin{a non-empty Zariski-open set, such that}\;\;\forall f\in\mathcal{O},\;\dim(\h\cdot f)=\frac{1}{2}\dim(\g\cdot f).
\end{equation} 
We will say that $f\in\g^{\ast}$ satisfies the lagrangian condition (w.r.t $\h$) if $\dim(\h\cdot f)=\frac{1}{2}\dim(\g\cdot f)$. Alternatively and depending on the context,  we will call $\h$ \textit{lagrangian} with respect to $f$ if it's simultaneously isotropic and coisotropic with respect to $B_f$, the Kostant form for $f$. Recall that if the $H\cdot f$ orbits are lagrangian in the $G\cdot f$ orbits (as Poisson submanifolds), then $(U_{\mathbb{C}}(\g)/U_{\mathbb{C}}(\g)\h_{-if})^{\h}$ is commutative (for this result see ~\cite{CG1} $\mathcal{x}$ 5, Theorem 5.4 and Corollary 5.5). Let $l\in\lambda+\h^{\bot}$ and $\g(l)$ be the stabilizer of $\l$. Such an $l$ will be called \textit{regular} if $\dim(\g(l))$ is minimal among the dimensions $\dim(\g(\xi)),\;\xi\in\lambda +\h^{\bot}$. The regular elements form a non-empty Zariski-open set. Furthermore, an $l\in\lambda+\h^{\bot}$ will be called \textit{generic} if it is regular and satisfies the lagrangian condition. The set of generic elements  is still a Zariski-open set. We suppose here that this set is non-empty. For the proof of the next result see ~\cite{FUJI2} Theorem 1, pp.757.
\newtheorem{ok}{Theorem}[section]
\begin{ok}[~\cite{FUJI3}]\label{fujisimple}
Let $\g$ be a finite dimensional Lie algebra, $\h\subset\g$ a subalgebra and $\lambda$ a character of $\h$. Suppose that generically  $\tau_{\lambda}=Ind(G\uparrow H,\chi_{\lambda})$ is of finite multiplicities. Then for generic $l\in \lambda+\h^{\bot}$ and for $A\in(U_{\mathbb{C}}(\g)/U_{\mathbb{C}}(\g)\h_{il})^{\h}$, $\mathrm{d}\tau_l(\overline{A})(\alpha_l)$ is a multiple of $\alpha_l$, and there is defined a character $\lambda_l:\;(U_{\mathbb{C}}(\g)/U_{\mathbb{C}}(\g)\h_{il})^{\h}\longrightarrow \mathbb{C}$ such that
\begin{equation}\label{rep character}
\mathrm{d}\tau_l(\overline{A})(\alpha_l)=\overline{\lambda_l(A)}\alpha_l.
\end{equation}
\end{ok}

\textbf{Deformation quantization (\cite{K}).} If $X$ is a Poisson manifold, a star-product on $C^{\infty}(X)$ is an $\mathbb{R}[[\epsilon]]$-bilinear map $C^{\infty}(X)[[\epsilon]]\times C^{\infty}(X)[[\epsilon]]\longrightarrow C^{\infty}(X)[[\epsilon]],\; (f,g)\mapsto f\ast g$, such that for $f,g,h \;\; \in C^{\infty}(X)$, $f\ast g= f\cdot g + \sum_{i=1}^{\infty}B_i(f,g)\epsilon^i$, the $B_i$ being bidifferential operators, $(f\ast g)\ast h = f\ast(g\ast h)$ and $f\ast 1 = 1\ast f =f$. To write his $\ast-$ product formula, Kontsevich used a family $\mathbf{Q}_{n,2}$ of graphs $\Gamma$ described as follows: The set $V(\Gamma)$ of vertices of $\Gamma$ is the disjoint union of two ordered sets $V_1(\Gamma)$ and $V_2(\Gamma)$, isomorphic to $\{1,\ldots,n\}$ and $\{1,2\}$ respectively. Their elements are called \textsl{type I} vertices, for $V_1(\Gamma)$, and \textsl{type II} vertices, for $V_2(\Gamma)$. The set $E(\Gamma)$ of edges of the graph is finite and all elements of $E(\Gamma)$ are oriented. If $S(r)$ is the set of edges leaving from $r\in V_1(\Gamma)$ then $\#S(r)=2$ and $\sum_{r\in V_1(\Gamma)\cap V_2(\Gamma)}S(r)=2n$, that is no edge leaves from a type II vertex. Furthermore every set $S(r)$ is ordered and no loops or double edges are allowed. $E(\Gamma)$ is ordered in a compatible way with the order in $V_1(\Gamma)$, and $S(r),\;r\in V_1(\Gamma)$. One associates a differential operator to each such graph: If $\{x_1,\ldots,x_k\}$ are coordinate variables of $X$ and $\pi$ is the Poisson bivector describing the Poisson structure on $X$, then to $\pi$ and a graph $\Gamma\in \mathbf{Q}_{n,2}$ we associate a differential operator $B_{\Gamma}^{\pi}$ in the following way: Let $L:\;E(\Gamma)\longrightarrow [1,k]:=\{1,\ldots,k\}$ be a labelling function of the edges of the graph. Fix a vertex $r\in [1,n]$ and let $S(r)=\{e_r^1,\ldots,e_r^2\}$ be the ordered set of edges leaving $r$. Associate the function $\pi^{L(e_r^1),L(e_r^2)}$ to $r$, where $\pi^{ij}$ is the matrix describing the Poisson structure. On each vertex of $V_2(\Gamma)$ we associate respectively a function $F,G\in C^{\infty}(X)$. To the $p^{th}-$ edge of $S(r)$, associate the partial derivative w.r.t the coordinate variable $x_{L(e_r^p)}$.This derivative is acting on the  function associated to the vertex $w\in V_1(\Gamma)\cup V_2(\Gamma)$ where the edge $e_r^p$ arrives. Let $(p,m)\in E(\Gamma)$ represent an oriented edge from $p$ to $m$. Then if $\{\overline{1},\overline{2}\}=V_2(\Gamma)$,
\begin{equation}\label{bidiff op}
B_{\Gamma}^{\pi}(F,G):= \sum_{L:E(\Gamma)\rightarrow [1,k]}\left[\prod_{r=1}^{\#(V_1(\Gamma))}\left(\prod_{\delta\in E(\Gamma),\;\delta=(\cdot,r)}\partial_{L(\delta)}\right)\pi^{L(\upsilon_r^1)L(\upsilon_r^2)}\right]\times
\end{equation}
\[\times\left(\prod_{\delta\in E(\Gamma),\;\delta=(\cdot,\overline{1})}\partial_{L(\delta)}\right)(F)\times\left(\prod_{\delta\in E(\Gamma),\;\delta=(\cdot,\overline{2})}\partial_{L(\delta)}\right)(G).\]
The last ingredient for Kontsevich's $\ast-$ product formula, the coefficient $\omega_{\Gamma}$. Let $\mathcal{H}=\{z\in \mathbb{C}/ \mathfrak{Im}(z)\geq 0\}$ be the upper-half plane and let $\mathcal{H}^+=\{z\in \mathbb{C}/ \mathfrak{Im}(z)> 0\}$. Embed an admissible graph $\Gamma$ in $\mathcal{H}$ by putting the type II vertices on the real axis and letting the type I vertices move in $\mathcal{H}^+$. Let $\widehat{C}_{n,\overline{m}}^+$ be the connected configuration manifold  $\widehat{C}_{n,\overline{2}}^+:=\{(z_1,\ldots,z_n,z_{\overline{1}},z_{\overline{2}})\in\widehat{C}_{n,\overline{2}}/\forall i<j,\;\; z_{\overline{1}}<z_{\overline{2}}\}$. Consider now the manifold $\widehat{C}_{2,0}$ (the hat denotes its invariance under horizontal translations and dilations), and a map, called the \textsl{angle map}, $\phi: \widehat{C}_{2,0}\longrightarrow \mathbb{R}/2\pi\mathbb{Z},\;
(z_1,z_2)\mapsto angle(<z_1,+i\infty>,<z_1,z_2>),$ where $<z_1,+i\infty>$ stands for the geodesic (with respect to the Poincare metric) passing by $z_1$ and $+i\infty$, and $<z_1,z_2>$ stands for the geodesic passing by $z_1$ and $z_2$. For $e=(z_i,z_j)$ an edge of $\Gamma$, let $p_{e}: \widehat{C}_{n,\overline{m}}\longrightarrow \widehat{C}_{2,0},\;(z_1,\ldots,z_n,z_{\overline{1}},\ldots,z_{\overline{m}})\mapsto (z_i,z_j)$ be the natural projection. One then defines a form on $\widehat{C}_{n,\overline{m}}$, setting $p^{\ast}_{e}(\mathrm{d}\phi)=:\mathrm{d}\phi_{e} \in \Omega^1(\widehat{C}_{n,\overline{m}})$. Then define $\Omega_{\Gamma}$ to be the form $\Omega_{\Gamma}:=\bigwedge_{e\in E(\Gamma)}\mathrm{d}\phi_e$ and finally the Kontsevich coefficient  $\omega_{\Gamma}:=\frac{1}{(2\pi)^{2n+\overline{m}-2}}\int_{\widehat{C}_{n,\overline{m}}^+}\Omega_{\Gamma}$. 
\newtheorem{okta}[ok]{Theorem}
\begin{okta}[~\cite{K}]\label{kon}
Let $X,\pi$ be a Poisson manifold. Then for $f,g\in C^{\infty}(X)$, the operator $f\ast_Kg:=fg+\sum_{n=1}^{\infty}\epsilon^n\left(\frac{1}{n!}\sum_{\Gamma\in\mathbf{Q_{n,2}}}\omega_{\Gamma}B_{\Gamma}^{\pi}(f,g)\right)$ is an associative product. 
\end{okta}
\textbf{Biquantization, colored graphs, bidifferential operators and coefficients (\cite{CF2},\cite{CF3},\cite{CT}).} In two fundamental papers ~\cite{CF2},~\cite{CF3}, A. Cattaneo and G. Felder, expanded the formality Theorem of M. Kontsevich for the case of one coisotropic submanifold $C$ of a Poisson manifold $X$. We adjust their results in our setting in the sense of \cite{CT}. Let $(X,\pi)$ be a Poisson manifold and $C$ a submanifold. Then $C$ is called \textsl{coisotropic} if the ideal $I(C)\subset C^{\infty}(X)$ of functions vanishing on $C$ is a Poisson ideal of $C^{\infty}(X)$. It is clear that in the Lie case, for a subalgebra $\h\subset\g$, the annihilator $\h^{\bot}$ (and also $\g^{\ast}$) is a coisotropic submanifold of $\g^{\ast}$ with the natural Poisson structure induced by the Lie bracket, $\pi=\frac{1}{2}[\cdot,\cdot]$.
Let now the edges of a graph be colored, either as $(+)$ or as $(-)$. Double edges are not allowed, meaning edges with the same color, source and target. Let $\mathfrak{q}$ be a supplementary space of $\h$, that is $\g=\h\oplus\mathfrak{q}$. Let also $\{H_1,H_2,\ldots, H_t\}$ be a basis for $\h$ and $\{Q_1,\ldots, Q_r\}$ a basis for $\mathfrak{q}$. We identify spaces $\mathfrak{q}^{\ast}\simeq\g^{\ast}/\h^{\ast}\simeq\h^{\bot}$. Graphically, the color $(-)$ will be represented with a dotted edge and the color $(+)$ will be represented with a straight edge. Recall that in order to construct the differential operator $B_{\Gamma}$ corresponding to a given graph, we associated a coordinate variable to each edge of the graph. 
In terms of $\g^{\ast}$, let $\{H_1^{\ast},H_2^{\ast},\ldots, H_t^{\ast}\}$ be the corresponding basis of $\h^{\ast}$ and $\{Q_1^{\ast},\ldots, Q_r^{\ast}\}$ the basis of $\mathfrak{q}^{\ast}$. Let $(x_i^{\ast})_{i=1,\ldots,n}$ be the coordinates relatively to the basis $\{H_1^{\ast},\ldots,H_t^{\ast},Q_1^{\ast},\ldots,Q_r^{\ast}\}$ and let $\Gamma$ be an admissible graph with two colors. To $\Gamma$ is associated a 2-colored 1-form $\overline{\Omega}_{\Gamma}$ and a 2-colored coefficient $\overline{\omega}_{\Gamma}$ setting $\phi_+:\;C_{2,0}\longrightarrow \mathbb{R}/2\pi\mathbb{Z}$ to be the function 
$\phi_+(z_1,z_2):=arg(z_1-z_2) +arg(z_1-\overline{z_2})$ and $\phi_-:\;C_{2,0}\longrightarrow \mathbb{R}/2\pi\mathbb{Z}$ be the function $\phi_-(z_1,z_2):=arg(z_1-z_2)-arg(z_1-\overline{z_2}).$ These functions extend to $\widehat{C}_{2,0}$ and thus can be used for definitions analogous to Kontsevich's. Note that $\mathrm{d}\phi_+(z_1,z_2)=\mathrm{d}\phi_-(z_2,z_1).$
This means that for the calculation of $\overline{\omega}_{\Gamma}$, an edge is equivalent with the reversed edge carrying the opposite color. The form $\overline{\Omega}_{\Gamma}$ of a 2-colored graph $\Gamma$ is similarly defined as $\overline{\Omega}_{\Gamma}:=\wedge_{e\in E(\Gamma)}\mathrm{d}\phi_{\cdot,e}$ where $\mathrm{d}\phi_{+,e}=p_e^{\ast}(\mathrm{d}\phi_+),\;\mathrm{d}\phi_{-,e}=p_e^{\ast}(\mathrm{d}\phi_-)$, when $e\in E(\Gamma)$ has color $+/-$ respectively, and the colored coefficient is  $\overline{\omega}_{\Gamma}:=\frac{1}{(2\pi)^{2n}}\int_{\widehat{C}_{n,2}^+}\overline{\Omega}_{\Gamma}$. If $\mathbf{Q}_{n,2}^{(2)}$ denotes the set of admissible graphs with two colors and two type II vertices, the formula (\ref{bidiff op} of $B(\Gamma)$ in this case has to be modified in the following sense: For $e\in E(\Gamma)$, let $c_e\in\{+,-\}$ be its color and let $L:\;E(\Gamma)\longrightarrow \{1,\ldots,t,t+1,\ldots t+r\}$, where $t=\dim(\h^{\ast}),\;r=\dim(\mathfrak{q}^{\ast})$, satisfying $L(e)\in\{1,\ldots,t\}\;\;\textlatin{if}\;c_e=-\;\;,\;\;L(e)\in\{t+1,\ldots,t+r\}\;\;\textlatin{if}\;c_e=+$ be a 2-colored labelling function. By its form, variables of $\h^{\ast}$ are associated to the color $(-)$ and variables of $\mathfrak{q}^{\ast}$ are associated to $(+)$. Then, for $F,G\in S(\g)$, the formula of \textsl{colored} bidifferential operators $B_{\Gamma}^{\pi}$ is the same as in (\ref{bidiff op}) but using the 2-colored labelling function $L$ that we just described and the family $\mathbf{Q}_{n,2}^{(2)}$ instead of $\mathbf{Q}_{n,2}$. Theorem \ref{kon} can be generalized in the following sense: The product $\ast_{CF,\epsilon}:\;C^{\infty}(C)[\epsilon]\times C^{\infty}(C)[\epsilon]\longrightarrow C^{\infty}(C)[\epsilon]$ given by the formula $F\ast_{CF,\epsilon} G:=F\cdot G+\sum_{n=1}^{\infty}\epsilon^n\left(\frac{1}{n!}\sum_{\Gamma\in\mathbf{Q^{(2)}_{n,2}}}\overline{\omega}_{\Gamma}B_{\Gamma}^{\pi}(F,G)\right)$, is associative and will be called the \textsl{Cattaneo-Felder} product. We now specify some special colored graphs that we will use (see ~\cite{CT} $\mathcal{x}$ 1.3, 1.6 and ~\cite{BAT} $\mathcal{x}$ 2.3). They are colored graphs that have an edge with no end colored as $(-)$ (associated to variables in $\mathfrak{q}^{\ast}$). We'll say that this edge "points to $\infty$", denote it by $e_{\infty}$. Denote this family of graphs as $\mathbf{Q}^{\infty}_{n,\overline{m}}$. 

 \textbf{Bernoulli.} The Bernoulli type graphs with $i$ type I vertices, $i\in\mathbb{N}$, will be denoted by $\mathcal{B}_i$. They derive the function $F$ $i$ times, have $2i$ edges and leave an edge towards $\infty$. These conditions imply the existence of a vertex $s\in V_1(\Gamma)$ that receives no edge, called the \textsl{root} of $\Gamma$.

 \textbf{Wheels.} The wheel type graphs with $i$ type I vertices, $i\in\mathbb{N}$, will be denoted by $\mathcal{W}_i$. They derive the function $F$ $i$ times, have $2i$ edges and leave no edge to $\infty$. 

 \textbf{Bernoulli attached to a wheel.} Graphs of this type with $i$ type I vertices, $i\in\mathbb{N}$, will be denoted by $\mathcal{BW}_i$. They derive the function $F$ $i-1$ times and leave an edge to $\infty$.
For an $\mathcal{W}_m-$ type graph $W_m$ attached to a $\mathcal{B}_l-$ type graph $B_l$, we'll write $B_lW_m\in\mathcal{B}_l\mathcal{W}_m$. Obviously $\mathcal{B}_l\mathcal{W}_m\subset\mathcal{BW}_{l+m}$.

Let us now give the definition of the reduction algebra without character. Let $\{e^1_l,\ldots,e^k_l\}$ be the ordered set of edges leaving the vertex $l\in V_1(\Gamma)$ of a colored graph $\Gamma_{\infty} \in\mathbf{Q}^{\infty}_{s,1}$. For such a $\Gamma$ and using the notation $H_i^{\ast}:=\partial_i$, let $B_{\Gamma}:\;S(\mathfrak{q})\longrightarrow S(\mathfrak{q})\otimes \h^{\ast},\;F\mapsto B_{\Gamma}(F)$ be the operator defined by the formula
\begin{equation}
B_{\Gamma}(F)= \sum_{\substack{L:E(\Gamma)\rightarrow [[1,t+r]]\\ \textlatin{L colored}\\}}\left[\prod_{r=1}^{n}\left(\prod_{\in E(\Gamma),\;e=(\cdot,r)}\partial_{L(e)}\right)\pi^{L(e_r^1)L(e_r^2)}\right]\times\left( \prod_{\substack{e\in E(\Gamma)\\ e=(\cdot,\overline{1})\\}}\partial_{L(e)}F\right)\otimes H^{\ast}_{L(e_{\infty})}
\end{equation}
 We denote as $d^{(\epsilon)}_{\h^{\bot},\mathfrak{q}}:\;S(\mathfrak{q})[\epsilon]\longrightarrow S(\mathfrak{q})[\epsilon]\otimes\h^{\ast}$ the differential operator $d^{(\epsilon)}_{\h^{\bot},\mathfrak{q}}=\sum_{i=1}^{\infty}\epsilon^i d^{(i)}_{\h^{\bot},\mathfrak{q}}$ where $d^{(i)}_{\h^{\bot},\mathfrak{q}}=\sum_{\Gamma\in \mathcal{B}_i\cup \mathcal{BW}_i}\overline{\omega}_{\Gamma}B_{\Gamma}$ and define the reduction algebra $H^0_{(\epsilon)}(\h^{\bot},d^{(\epsilon)}_{\h^{\bot},\mathfrak{q}})$ of polynomials in $\epsilon$ as the vector space of solutions of the equation 
\begin{equation}\label{reduction equations}
d^{(\epsilon)}_{\h^{\bot},\mathfrak{q}}(F_{(\epsilon)})=0,\;F_{(\epsilon)}\in S(\mathfrak{q})[\epsilon]
\end{equation}
 equipped with the $\ast_{CF,\epsilon}-$ product, (which is associative on $H^0_{(\epsilon)}(\h^{\bot},d^{(\epsilon)}_{\h^{\bot},\mathfrak{q}})$ by ~\cite{CF3}).

$\mathbf{H^0(\h^{\bot},d_{\h^{\bot},\mathfrak{q}}).}$ The system (\ref{reduction equations}) can also be written as homogeneous equations this time not using the degree on $\epsilon$, but instead the degree $\mathrm{deg}_{\mathfrak{q}}$ of each operator $B_{\Gamma}$. This results in the same system of homogeneous equations for a function $F=\sum_{i=0}^nF^{(n-i)}$ in $S(\mathfrak{q})$ with each $F^{(k)}$ being a homogeneous polynomial of degree k ($\mathrm{deg}_{\mathfrak{q}}(F^{(k)})=k$). For $d_{\h^{\bot},\mathfrak{q}}:\;S(\mathfrak{q})\longrightarrow S(\mathfrak{q})\otimes \h$, the differential operator $d_{\h^{\bot},\mathfrak{q}}= \sum_{i=1}^{\infty}d^{(i)}_{\h^{\bot},\mathfrak{q}}$ where $d^{(i)}_{\h^{\bot},\mathfrak{q}}=\sum_{\Gamma\in \mathcal{B}_i\cup \mathcal{BW}_i}\overline{\omega}_{\Gamma}B_{\Gamma}$, we denote as $H^0(\h^{\bot},d_{\h^{\bot},\mathfrak{q}})$ the vector space of polynomials $F$, solutions of the equation $d_{\h^{\bot},\mathfrak{q}}(F)=0$. The Cattaneo-Felder construction without $\epsilon$ is still valid for polynomial functions and defines an associative product $\ast_{CF}$ on $H^0(\h^{\bot},d_{\h^{\bot},\mathfrak{q}})$. In the sequence $H^0(\h^{\bot},d_{\h^{\bot},\mathfrak{q}})$ will stand for the algebra $\left(H^0(\h^{\bot},d_{\h^{\bot},\mathfrak{q}}),\ast_{CF}\right)$.

$\mathbf{H^0_{(\epsilon)}(\h_{\lambda}^{\bot},d^{(\epsilon)}_{\h^{\bot}_{\lambda},\mathfrak{q}})}$ \textbf{and} $\mathbf{H^0_{(\epsilon=1)}(\h_{\lambda}^{\bot},d^{(\epsilon=1)}_{\h^{\bot}_{\lambda},\mathfrak{q}})}$. Let $\h_{\lambda}^{\bot}:=\{f\in\g^{\ast}/f|_{\h}=-\lambda\}=-\lambda+\h^{\bot}\subset\g^{\ast}$. On $-\lambda+\h^{\bot}$ define the reduction algebra $H^0_{(\epsilon)}(\h_{\lambda}^{\bot},d^{(\epsilon)}_{\h^{\bot}_{\lambda},\mathfrak{q}})$ using the deformation parameter $\epsilon$ with the same way as in the vector space case. We denote by $H^0_{(\epsilon)}(\h_{\lambda}^{\bot},d^{(\epsilon)}_{\h^{\bot}_{\lambda},\mathfrak{q}})$ the algebra of non $\mathrm{deg}_{\mathfrak{q}}-$ homogeneous $\epsilon-$ polynomials $P_{(\epsilon)}$, solutions of the equation $d^{(\epsilon)}_{\h^{\bot}_{\lambda},\mathfrak{q}}(P_{(\epsilon)})=0$, equipped with the $\ast_{CF,\epsilon}$ product, and the specialized algebra of $\h_{\lambda}^{\bot}$ as $H^0_{(\epsilon=1)}(\h_{\lambda}^{\bot},d^{(\epsilon=1)}_{\h_{\lambda}^{\bot},\mathfrak{q}}):=\left(H^0_{(\epsilon)}(\h_{\lambda}^{\bot},d^{(\epsilon)}_{\h_{\lambda}^{\bot},\mathfrak{q}})/<\epsilon-1>\right)$. The Cattaneo-Felder product on $H^0_{(\epsilon=1)}(\h_{\lambda}^{\bot},d^{(\epsilon=1)}_{\h_{\lambda}^{\bot},\mathfrak{q}})$ will be denoted also as $\ast_{CF,(\epsilon=1)}$. In the same way we will consider the associative algebra$\left(H^0(\h_{\lambda}^{\bot},d_{\h_{\lambda}^{\bot},\mathfrak{q}}),\ast_{CF}\right)$ denoted simply as $H^0(\h_{\lambda}^{\bot},d_{\h_{\lambda}^{\bot},\mathfrak{q}})$.

An important feature introduced in \cite{CF2}-\cite{CF3} is the biquantization (and triquantization) diagrams. In this paper we are using the biquantization diagrams of $\g^{\ast}, \h^{\bot}_{\lambda}$, of $-f+\mathfrak{b}^{\bot}, \h^{\bot}_{\lambda}$ and the triquantization diagram of $\g^{\ast}, \h_{\lambda}^{\bot}, -f+\mathfrak{b}^{\bot}$, for $\mathfrak{b}$ polarization with respect to $f$. More specifically there are bimodule structures of the reduction space (corresponding to the intersection of the two coisotropic submanifolds) at the corner of each such diagram with respect to the reduction algebras in the axes. For example in the biquantization diagram of $-f+\mathfrak{b}^{\bot}, -\lambda+\h^{\bot}$, there is a right $H^0_{(\epsilon)}(\h_{\lambda}^{\bot},d^{(\epsilon)}_{\h^{\bot}_{\lambda},\mathfrak{q}})-$ module and a left $ H^0_{(\epsilon)}(-f + \mathfrak{b}^{\bot},d^{(\epsilon)}_{\mathfrak{b}^{\bot}})-$ module structure of $H^0_{(\epsilon)}(-f+(\h+\mathfrak{b})^{\bot},d^{(\epsilon)}_{\h^{\bot},\mathfrak{b}^{\bot}})$ (reduction space at the origin of this diagram).

The paper is structured as follows: Section 2.1 contains Theorem \ref{ctb} that constructs a character for the associative algebra $(U_{(\epsilon)}(\g)/U_{(\epsilon)}(\g)\h_{f})^{\h}$. We explain how the lagrangian condition permits us to construct many characters for this algebra. Section 2.2 is introducing the algebra $\mathbb{S}(\g,h,\lambda)$ of \textsl{affine symbols}. It is equipped with a more or less natural Poisson structure and we use it to relate the commutativity of this algebra with the commutativity of $H^0_{(\epsilon)}(\h_{\lambda}^{\bot},d^{(\epsilon)}_{\h^{\bot}_{\lambda},\mathfrak{q}})$. In section 2.3 we prove Theorem \ref{fujireal} which is the corresponding result of Theorem \ref{fujisimple} in the real case for the specialization algebra of differential operators $\mathbb{D}(\,\h,\lambda)$, and prepare the ground for the proof of the central result of this paper, Theorem \ref{central}. It states that the characters constructed in Theorems \ref{ctb} and \ref{fujireal} actually coincide. In section 3, a detailed example of a specific nilpotent Lie algebra, originally treated by Fujiwara, is used to compute the same characters but via deformation quantization. The benefit is that an explicit formula for characters of  $(U_{\mathbb{C}}(\g)/U_{\mathbb{C}}(\g)\h_{i\lambda})^{\h}$ is obtained thanks to these techniques. An interesting and unknown term with rational coefficients appears in this formula.

\section{Construction and comparison of characters.}

\subsection{Construction.}
Let $\mathfrak{b}$ be a polarization of $f\in\g^{\ast}$ and $\h\subset\g$ a subalgebra .We'll call $\mathfrak{q}$ a transversal supplementary of $\h$ with respect to $\mathfrak{b}$ if it satisfies $\mathfrak{b}=\mathfrak{b}\cap\h\oplus\mathfrak{b}\cap\mathfrak{q}$ and of course $\g=\h\oplus\mathfrak{q}$. It will be denoted as $\mathfrak{q}_{\mathfrak{b}}$. We'll need the following from ~\cite{CT}.

\newtheorem{p}[ok]{Proposition}
\begin{p}[~\cite{CT}]\label{ct}
Let $\g$ be a real Lie algebra, $\h\subset\g$ a subalgebra and $f\in\g^{\ast}$ s.t $f([\h,\h])=0$. Suppose that there is $\mathfrak{b}$ a polarization with respect to $f$ and that $\h$ is lagrangian with respect to~$f$. Construct the biquantization diagram of $-f+\h^{\bot},\;-f+\mathfrak{b}^{\bot}$, with $ -f+\mathfrak{b}^{\bot}$ at the vertical axis. Let $\mathfrak{q}_{\mathfrak{b}}$ be a transversal supplementary space of $\h$ in $\g$ and let $P\in H^0_{(\epsilon)}(\h_f^{\bot},d^{(\epsilon)}_{\h_f^{\bot},\mathfrak{q}_{\mathfrak{b}}})$.  Let also 
\[\gamma:  H^0_{(\epsilon)}(\h_f^{\bot},d^{(\epsilon)}_{\h_f^{\bot},\mathfrak{q}_{\mathfrak{b}}})\longrightarrow H^0_{(\epsilon)}(-f+(\h+\mathfrak{b})^{\bot},d^{(\epsilon)}_{\h^{\bot},\mathfrak{b}^{\bot}}),\]
\[P\mapsto P\ast_{1,\mathfrak{b}}1\]
where $\ast_{1,\mathfrak{b}}$ is the left $H^0_{(\epsilon)}(\h_{f}^{\bot},d^{(\epsilon)}_{\h_{f}^{\bot},\mathfrak{q}_{\mathfrak{b}}})-$ module structure of Cattaneo-Felder for the reduction space $H^0_{(\epsilon)}(-f+(\h+\mathfrak{b})^{\bot},d^{(\epsilon)}_{\h^{\bot},\mathfrak{b}^{\bot}})$.
Then $\gamma$ is a character for $H^0_{(\epsilon)}(\h_f^{\bot},d^{(\epsilon)}_{\h_f^{\bot},\mathfrak{q}_{\mathfrak{b}}})$.
\end{p}
\textit{Proof}. Either in ~\cite{CT} or in ~\cite{BAT} Prop. 4.3. It is based in proving that $ H^0_{(\epsilon)}(-f + \mathfrak{b}^{\bot},d^{(\epsilon)}_{\mathfrak{b}^{\bot}})=\mathbb{R}[\epsilon]$ and that also $H^0_{(\epsilon)}(-f+(\h+\mathfrak{b})^{\bot},d^{(\epsilon)}_{\h^{\bot},\mathfrak{b}^{\bot}})=\mathbb{R}[\epsilon]$. $\diamond$

Recall that in \cite{CT} $\mathcal{x}$ 6.2.1 the authors gave a precise definition of an 8-colored 1-form $\theta_{j_1j_2j_3}$ which in an abstract triquantization diagram interpolates the two 4-colored 1-forms of biquantization at the corners of the diagram. In $\mathcal{x}$ 6.2.2 they used the form $\theta_{j_1j_2j_3}$ in a triquantization diagram of $f+\mathfrak{b}_1^{\bot},\h,f+\mathfrak{b}^{\bot}_2$, where $\mathfrak{b}_1,\mathfrak{b}_2$ are polarizations of $f$ in normal intersection position ($\h\cap(\mathfrak{b}_1+\mathfrak{b}_2)=\h\cap\mathfrak{b}_1+\h\cap\mathfrak{b}_2$), to prove that the character computed in Proposition \ref{ct} is independent of the chosen polarization. Let now
\[c_e=(+,+,+)\longrightarrow X_e^{\ast}\in (\g\cap\h\cap\mathfrak{b})^{\ast},\;\;c_e=(-,+,+)\longrightarrow X_e^{\ast}\in(\h\cap\mathfrak{b})^{\ast}/\g^{\ast},\]
 \begin{equation}\label{coloring with polarization 3}
c_e=(+,-,+)\longrightarrow X_e^{\ast}\in (\g\cap\mathfrak{b})^{\ast}/(\h\cap\mathfrak{b})^{\ast},\;c_e=(-,-,+)\longrightarrow X_e^{\ast}\in\mathfrak{b}^{\ast}/(\g+\h)^{\ast},
\end{equation}
\[c_e=(+,+,-)\longrightarrow X_e^{\ast}\in (\g\cap\h)^{\ast}/(\h\cap\mathfrak{b})^{\ast},\;\;c_e=(-,+,-)\longrightarrow X_e^{\ast}\in\h^{\ast}/(\g+\mathfrak{b})^{\ast},\]
\[c_e=(+,-,-)\longrightarrow X_e^{\ast}\in \g^{\ast}/(\h+\mathfrak{b})^{\ast},\;\;c_e=(-,-,-)\longrightarrow X_e^{\ast}\in\g^{\ast}/(\g+\h+\mathfrak{b})^{\ast}.\]
be a coloring of graphs in the triquantization diagram with $\g^{\ast}$ at the left vertical axis, $-f+\h^{\bot}$ at the horizontal axis and $-f+\mathfrak{b}^{\bot}$ at the right vertical axis. Obviously the colors in the right column of (\ref{coloring with polarization 3}) do not appear, but initially we need to consider ourselves in a triquantization setting. Let $T_1,T_2$ be the operators $T_1:\;H^0_{(\epsilon)}(\g^{\ast},d^{(\epsilon)}_{\g^{\ast}})\longrightarrow H^0_{(\epsilon)}(\h_{\lambda}^{\bot},d^{(\epsilon)}_{\g^{\ast},\h_{\lambda}^{\bot},\mathfrak{q}}),\;\; F \mapsto F\ast_1 1,$ and $T_2:\;H^0_{(\epsilon)}(\h^{\bot}_{\lambda},d^{(\epsilon)}_{\h^{\bot}_{\lambda},\mathfrak{q}}) \longrightarrow H^0_{(\epsilon)}(\h_{\lambda}^{\bot},d^{(\epsilon)}_{\g^{\ast},\h_{\lambda}^{\bot},\mathfrak{q}}),\;\;G\mapsto 1 \ast_2 G$, where $\ast_1$ corresponds to the left $H^0_{(\epsilon)}(\g^{\ast},d^{(\epsilon)}_{\g^{\ast}})-$ module structure of $H^0_{(\epsilon)}(\h_{\lambda}^{\bot},d^{(\epsilon)}_{\g^{\ast},\h_{\lambda}^{\bot},\mathfrak{q}})$ and $\ast_2$ corresponds to the right $H^0_{(\epsilon)}(\h^{\bot}_{\lambda},d^{(\epsilon)}_{\h^{\bot}_{\lambda},\mathfrak{q}}) -$ module structure of $H^0_{(\epsilon)}(\h_{\lambda}^{\bot},d^{(\epsilon)}_{\g^{\ast},\h_{\lambda}^{\bot},\mathfrak{q}})$ in the biquantization diagram of $\g^{\ast},\h^{\bot}_{\lambda}$. Let also $q(Y):= \det_{\g} \left(\frac{\sinh\frac{\mathrm{ad} Y}{2}}{\frac{adY}{2}}\right)$, and $\forall Y\in\g$ set $q_{(\epsilon)}(Y):=q(\epsilon X)$. Finally let $\beta:\;S(\g)\longrightarrow U(\g)$ denote the symmetrization, $\beta_{(\epsilon)}:\;S_{(\epsilon)}(\g)\longrightarrow U_{(\epsilon)}(\g)$ the deformed symmetrization (the PBW Theorem still holds for $S_{(\epsilon)}(\g),\;U_{(\epsilon)}(\g)$), and $\overline{\beta}_{\mathfrak{q},(\epsilon)}:\;S(\mathfrak{q})[\epsilon]\longrightarrow~(U_{(\epsilon)}(\g)/U_{(\epsilon)}(\g)\h_{\lambda}$ the quotient deformed symmetrization w.r.t a supplementary $\mathfrak{q}$. Recall that the main result of \cite{BAT1} was Theorem 5.1 stating that the map $\overline{\beta}_{\mathfrak{q},(\epsilon)}\circ\partial_{q_{(\epsilon)}^{\frac{1}{2}}}\circ \overline{T}_1^{-1}T_2:\;H^0_{(\epsilon)}(\h^{\bot}_{\lambda},d^{(\epsilon)}_{\h^{\bot}_{\lambda},\mathfrak{q}})\stackrel{\sim}{\longrightarrow}\left(U_{(\epsilon)}(\g)/U_{(\epsilon)}(\g)\h_{\lambda}\right)^{\h}$ is an non-canonical algebra isomorphism, where  $\overline{T}_1:=T_1|_{S(\mathfrak{q})[\epsilon]}$.
\newtheorem{okan}[ok]{Theorem}
\begin{okan}\label{ctb}
Let $\g$ be a real Lie algebra, $\h\subset\g$ a subalgebra, and $f\in\g^{\ast}$ s.t $\h$ is lagrangian with respect to $f$. Suppose also that there is a polarization $\mathfrak{b}$ of $f$. Let finally $\mathfrak{q}_{\mathfrak{b}}$ be a transversal supplementary space of $\h$.
The linear map of $\left(U_{(\epsilon)}(\g)/U_{(\epsilon)}(\g)\h_{f}\right)^{\h}$ to $\mathbb{R}[\epsilon]$ given by
\[\gamma_{CT}:\;\left(U_{(\epsilon)}(\g)/U_{(\epsilon)}(\g)\h_{f}\right)^{\h}\longrightarrow\mathbb{R}[\epsilon]\]
\begin{equation}\label{character formula for U}
u\mapsto \overline{T}_1^L\circ\bar{\beta}_{\mathfrak{q}_{\mathfrak{b}},(\epsilon)}^{-1}(u)(f)
\end{equation}
is a character i.e $\forall u,v\in (U_{(\epsilon)}(\g)/U_{(\epsilon)}(\g)\h_{f})^{\h},\;\;\gamma_{CT}(uv)=\gamma_{CT}(u)\gamma_{CT}(v)$.
\end{okan}
\textit{Proof.}
The proof is a combination of results of ~\cite{CT} with Prop. \ref{ct}.
Consider the triquantization diagram with $\g^{\ast}$ at the left vertical axis, $-f+\h^{\bot}$ at the horizontal axis and $-f+\mathfrak{b}^{\bot}$ at the right vertical axis.  Choose a transversal supplementary space $\mathfrak{q}_{\mathfrak{b}}$. Let $P\in H^0_{(\epsilon)}(\h_f^{\bot},d^{(\epsilon)}_{\h_f^{\bot},\mathfrak{q}_{\mathfrak{b}}})$ on the horizontal axis approach the right corner of the diagram. By $\cite{CT}$, and Prop. \ref{ct} we calculate there a character for $H^0_{(\epsilon)}(\h_{f}^{\bot},d^{(\epsilon)}_{\h_{f}^{\bot},\mathfrak{q}_{\mathfrak{b}}})$. This character is $\gamma:\;P\mapsto \gamma(P)$ where $\gamma(P)=P\ast_{1,\mathfrak{b}} 1$ and $\ast_{1,\mathfrak{b}}$ is the left $H^0_{(\epsilon)}(\h_{f}^{\bot},d^{(\epsilon)}_{\h_{f}^{\bot},\mathfrak{q}_{\mathfrak{b}}})-$ module structure of $H^0_{(\epsilon)}(-f+ (\h+\mathfrak{b})^{\bot},d^{(\epsilon)}_{\h^{\bot},\mathfrak{b}^{\bot}})$. We prefer to denote the character on the right corner as $P\mapsto T_2^R(P):=P\ast_{1,\mathfrak{b}} 1$. The action of exterior graphs deriving in the concentration is trivial, since we derive $T_2^R(P)$, a constant function on $-f+(\h+\mathfrak{b})^{\bot}$.
\begin{figure}[h!]
\begin{center}
\includegraphics[width=7cm]{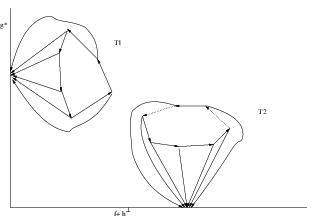}
\caption\footnotesize{The operators $\overline{T}_1^L$ and $T_2^L$ at the left corner of the triquantization diagram for $\g^{\ast},-f+\h^{\bot},-f+\mathfrak{b}^{\bot}$.}
\end{center}
\end{figure}
Move now $P$ at the left corner of the diagram. There we denote by $\ast_{2,\h}$ the right $H^0_{(\epsilon)}(\h_{f}^{\bot},d^{(\epsilon)}_{\h_{f}^{\bot},\mathfrak{q}_{\mathfrak{b}}})-$ module structure of $H^0_{(\epsilon)}(\h_{f}^{\bot},d^{(\epsilon)}_{\g^{\ast},\h_{f}^{\bot},\mathfrak{q}_{\mathfrak{b}}})$ and by $T_2^L(P)=1 \ast_{2,\h} P$ the operator composed of $\mathcal{W}$-type graphs on the left corner of the horizontal axis. Exterior graphs acting nontrivially on $T_2^L(P)$, are nonzero operators that have edges colored by $(+,+,-)$ and $(+,-,-)$ according to the coloring (\ref{coloring with polarization 3}). The 1-form $\omega^+_{k_1k_2k_3}(p,q)$ of 8 colors that we use to calculate the weight of such graphs, is zero if $q$ is in the corner. So only egdes of color $(+,-,-)$, i.e with variables in $(\h + \mathfrak{b})^{\bot}$  can derive nontrivially at this corner. Thus we have to restrict $T_2^L(P)$ in this direction and consider the restriction $T_2^L(P)|_{-f + (\h + \mathfrak{b})^{\bot}}$. Let $\Gamma_{int}$ be the family of possible interior graphs, i.e graphs appearing in $T_2^L$, and $\Gamma_{ext}$ the family of possible exterior graphs. Let also $A=\exp (\sum_{\Delta \in\Gamma_{ext}} \Delta)$ be the operator coming from the exterior contributions acting on the corner, and $B=\exp (\sum_{\Delta \in\Gamma_{int}}\Delta)$ be the operator corresponding to graphs in $T_2^L$. The operator $A$ consists of $\mathcal{W}$-type graphs because all the edges arriving at the corner have the same color. That is, the graphs in $A$ correspond to differential operators with constant coefficients. In this special case, we can write $A(B(P))|_{-f + (\h + \mathfrak{b})^{\bot}}=A(B(P)|_{-f + (\h + \mathfrak{b})^{\bot}})$ because $A$ derives in the same direction as the restriction. That is, the color $(+,-,-)$ that carry the edges of $A$ deriving $B$, corresponds to variables in $(\h+\mathfrak{b})^{\bot}$, the space where we have to restrict $A(B)$. Let $A(B(P)|_{-f + (\h + \mathfrak{b})^{\bot}})=c_f$ be the result of the evaluation of $A$ on $B(P)|_{-f + (\h + \mathfrak{b})^{\bot}}$.
By ~\cite{CT}, $c_f$ is a constant function on $-f+(\h+\mathfrak{b})^{\bot}$ as it coincides with the evaluation $T_2^R(P)(f)$ of the character $P\mapsto T_2^R(P)$ of  $H^0_{(\epsilon)}(\h_{f}^{\bot},d^{(\epsilon)}_{\h_{f}^{\bot},\mathfrak{q}_{\mathfrak{b}}})$ at the right corner of the diagram. Since $A$ is an invertible differential operator of constant coefficients, the fact that $c_f$ is a constant function means that $A$ acts trivially on $B(P)|_{-f + (\h + \mathfrak{b})^{\bot}}$. So at the left corner, the map
\[P\mapsto A(B(P)|_{-f + (\h + \mathfrak{b})^{\bot}})=B(P)|_{-f + (\h + \mathfrak{b})^{\bot}}=T_2^L(P)|_{-f+(\h+\mathfrak{b})^{\bot}}=(1\ast_{2,\h} P)|_{-f+(\h+\mathfrak{b})^{\bot}}\]
is a character of $H^0_{(\epsilon)}(\h_{f}^{\bot},d^{(\epsilon)}_{\h_{f}^{\bot},\mathfrak{q}_{\mathfrak{b}}})$. By Theorem 5.1 in \cite{BAT1}, we write $P$ as $P= (T_2^L)^{-1}\overline{T}_1^L\circ \bar{\beta}_{\mathfrak{q}_{\mathfrak{b}},(\epsilon)}^{-1}(u) \in H^0_{(\epsilon)}(\h_{f}^{\bot},d^{(\epsilon)}_{\h_{f}^{\bot},\mathfrak{q}_{\mathfrak{b}}})$ with $u\in \left(U_{(\epsilon)}(\g)/U_{(\epsilon)}(\g)\h_{f}\right)^{\h}$. Thus since the map $P\mapsto T_2^L(P)|_{-f+(\h+\mathfrak{b})^{\bot}}$ is a character for $H^0_{(\epsilon)}(\h_{f}^{\bot},d^{(\epsilon)}_{\h_{f}^{\bot},\mathfrak{q}_{\mathfrak{b}}})$, we get a character for $\left(U_{(\epsilon)}(\g)/U_{(\epsilon)}(\g)\h_{f}\right)^{\h}$, namely $\gamma_{CT}:\left(U_{(\epsilon)}(\g)/U_{(\epsilon)}(\g)\h_{f}\right)^{\h}\longrightarrow \mathbb{R}[\epsilon]$ with $\gamma_{CT}(u)= \overline{T}_1^L\circ\bar{\beta}_{\mathfrak{q}_{\mathfrak{b}},(\epsilon)}^{-1}(u)(f)$.  $\diamond$
\begin{figure}[h!]
\begin{center}
\includegraphics[width=9cm]{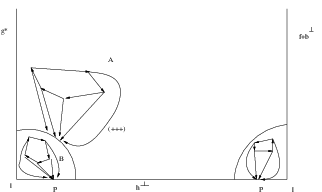}
\caption\footnotesize{The diagram constructing $\gamma_{CT}$.}
\end{center}
\end{figure}

Note for the future that the function  $(U_{(\epsilon)}(\g)/U_{(\epsilon)}(\g)\h_{f})^{\h}\ni u\mapsto \overline{T}_1^L\circ\overline{\beta}^{-1}_{\mathfrak{q}_{\mathfrak{b}},(\epsilon)}(u)$ is a constant function on $-f+(\h+\mathfrak{b})^{\bot}$. Also we didn't supposed so far that $(U_{(\epsilon)}(\g)/U_{(\epsilon)}(\g)\h_{f})^{\h}$ is commutative. From ~\cite{FLMM} the condition of existence of an $\mathcal{O}\subset \lambda+\h^{\bot}$ s.t $\forall l\in\mathcal{O}$ the $H\cdot l$ orbits are lagrangian, is equivalent to the condition $\mathbb{D}(\g,\h,\lambda)$ \textit{is commutative}. In view of this we can construct a family of characters.

\newtheorem{pronto}[ok]{Proposition}
\begin{pronto}\label{manychars}
Let $\g$ be a real nilpotent Lie algebra, $\h\subset\g$ a subalgebra and $G,H$ the corresponding connected and simply connected nilpotent Lie groups. Let also $\lambda\in\h^{\ast}$ a character of $\h$, $f\in\g^{\ast}$ s.t $f|_{\h}=\lambda$, $\mathfrak{b}$ a polarization of $f$ and $\mathfrak{q}_{\mathfrak{b}}$ a transversal supplementary space. Suppose that $\exists \mathcal{O}\subset \lambda+\h^{\bot}$ non-empty, Zariski-open subset such that $\forall l\in\mathcal{O}$, the $H\cdot l$ orbits are lagrangian. Then for every such $l$, there is a character of $(U_{(\epsilon)}(\g)/U_{(\epsilon)}(\g)\h_{\lambda})^{\h}$ namely $\gamma_{CT}^l:\left(U_{(\epsilon)}(\g)/U_{(\epsilon)}(\g)\h_{\lambda}\right)^{\h}\longrightarrow \mathbb{R}[\epsilon],\;\;u \mapsto \overline{T}_1^L\circ\bar{\beta}_{\mathfrak{q}_{\mathfrak{b}},(\epsilon)}^{-1}(u)(l)$.
\end{pronto}
\textit{Proof}. Since we are in the nilpotent case, there is always a polarization for the element $f$. Due to the lagrangian condition, Theorem \ref{ctb} gives us a character for each $l\in\mathcal{O}$. $\diamond$

\subsection{Affine symbols and applications.}
Let $\g$ be a real nilpotent Lie algebra, $\h\subset\g$ a subalgebra, $\lambda$ a character of $\h$ and $f\in\g^{\ast}$ s.t $f|_{\h}=\lambda$. Let $P_{(\epsilon)}\in H^0_{(\epsilon)}(\h^{\bot}_{\lambda},d^{(\epsilon)}_{\h^{\bot}_{\lambda},\mathfrak{q}})$ with $P_{(\epsilon)}=\sum_{i=0}^n\epsilon^iP_i$. Recall that the terms $P_i$ are not homogeneous and that $H^0_{(\epsilon)}(\h^{\bot}_{\lambda},d^{(\epsilon)}_{\h^{\bot}_{\lambda},\mathfrak{q}})\subset S(-f+\h^{\bot})\otimes\mathbb{R}[\epsilon]$. We introduce the notion of "affine symbol", that is a map $\sigma_{aff}:\;H^0_{(\epsilon)}(\h^{\bot}_{\lambda},d^{(\epsilon)}_{\h^{\bot}_{\lambda},\mathfrak{q}}) \longrightarrow \mathbb{R}[-f+\h^{\bot}]^{\h},$ with $P_{(\epsilon)}\mapsto P_0$.
\newtheorem{s}[ok]{Lemma}
\begin{s}
The symbol map $\sigma_{aff}$ takes values in $\mathbb{R}[-f+\h^{\bot}]=\left(S(\g)/S(\g)\h_{\lambda}\right)^{\h}$ and is an algebra map.
\end{s}
\textit{Proof.} Let $P_{(\epsilon)}\in H^0_{(\epsilon)}(\h^{\bot}_{\lambda},d^{(\epsilon)}_{\h^{\bot}_{\lambda},\mathfrak{q}})$.  Recall from the proof of Theorem 5.1 in \cite{BAT1}, that the first equation of the system 
$\sum_{\alpha}\left(\sum_{\Gamma^{\alpha}_{int},\Gamma^{\alpha}_{ext}}\left(\sum_{l,k,m}\left(B^m_{\Gamma^{\alpha}_{ext}}(B^k_{\Gamma^{\alpha}_{int}}(P_{n-l}))\right)\epsilon^{m+k+l}\right)\right)=0$, equivalent to a Stokes equation, was  
$B^1_{\Gamma^{\alpha}_{int}}(P_0)=~0$ which implies that $P_0\in \mathbb{R}[-f+\h^{\bot}]^{\h}$. Furthermore, let $\alpha=\sum_{k=0}^{\infty}\epsilon^k\alpha_k \in \mathbb{R}[\epsilon]$. Then $\alpha\cdot P_{(\epsilon)}=\sum_{k=0}^{\infty}\epsilon^k(\sum_{i+j=k}\alpha_iP_j)$ and $\alpha P_{(\epsilon)}\in H^0_{(\epsilon)}(\h^{\bot}_{\lambda},d^{(\epsilon)}_{\h^{\bot}_{\lambda},\mathfrak{q}})$ since $P_{(\epsilon)}\in H^0_{(\epsilon)}(\h^{\bot}_{\lambda},d^{(\epsilon)}_{\h^{\bot}_{\lambda},\mathfrak{q}})$. In this case, $\sigma_{aff}(\alpha P_{(\epsilon)})=\alpha_0 P_0$, and similarly $\sigma_{aff}(F_{(\epsilon)}+G_{(\epsilon)})=F_0+G_0$. For the product of $F_{(\epsilon)},G_{(\epsilon)}\in H^0_{(\epsilon)}(\h^{\bot}_{\lambda},d^{(\epsilon)}_{\h^{\bot}_{\lambda},\mathfrak{q}})$, let $F_{(\epsilon)}=\sum_{i=0}^n\epsilon^iF_i,\;G_{(\epsilon)}=\sum_{j=0}^p\epsilon^jG_j$. Then $\sigma_{aff}(F_{(\epsilon)}\ast_{CF}G_{(\epsilon)})=F_0G_0.\;\diamond$

Note however that $\sigma_{aff}$ is not an injective map, since $\forall F\in H^0_{(\epsilon)}(\h^{\bot}_{\lambda},d^{(\epsilon)}_{\h^{\bot}_{\lambda},\mathfrak{q}}),\;F\neq 0\Rightarrow \sigma_{aff}(\epsilon F)=~0$. Let now $\mathbb{S}(\g,\h,\lambda)$ be the algebra $\mathbb{S}(\g,\h,\lambda):=\{\sigma_{aff}(P_{(\epsilon)})/\;P_{(\epsilon)}\in H^0_{(\epsilon)}(\h^{\bot}_{\lambda},d^{(\epsilon)}_{\h^{\bot}_{\lambda},\mathfrak{q}})\}$.
It can be equipped with a Poisson structure: Consider two functions $P_{(\epsilon)}=\sum_{i=0}^n\epsilon^iP_i,\;Q_{(\epsilon)}=\sum_{j=0}^p\epsilon^jQ_j$, with $P_{(\epsilon)},Q_{(\epsilon)}\in H^0_{(\epsilon)}(\h^{\bot}_{\lambda},d^{(\epsilon)}_{\h^{\bot}_{\lambda},\mathfrak{q}})$. Then define 
\[\{P_0,Q_0\}:=\sigma_{aff}\left(\frac{P_{(\epsilon)}\ast_{CF} Q_{(\epsilon)}- Q_{(\epsilon)}\ast_{CF} P_{(\epsilon)}}{\epsilon}\right),\]
where $\sigma_{aff}(P_{(\epsilon)})=P_0,\sigma_{aff}(Q_{(\epsilon)})=Q_0$. This Poisson structure is well defined, that is, it doesn't depend on the choice of $P_{(\epsilon)} ,Q_{(\epsilon)}$. Indeed, let $\tilde{P}_{(\epsilon)}$ be another choice such that $\tilde{P}_{(\epsilon)}-P_{(\epsilon)}=\epsilon\cdot R$. Then $\tilde{P}_{(\epsilon)}\ast_{CF} Q_{(\epsilon)}- Q_{(\epsilon)}\ast_{CF} \tilde{P}_{(\epsilon)}=P_{(\epsilon)}\ast_{CF} Q_{(\epsilon)}- Q_{(\epsilon)}\ast_{CF} P_{(\epsilon)}+\epsilon [R, Q_{(\epsilon)}]$ with $ \epsilon [R, Q_{(\epsilon)}]\in O(\epsilon^2)$. Thus $\frac{[\tilde{P}_{(\epsilon)},Q_{(\epsilon)}]}{\epsilon}=\frac{[P_{(\epsilon)},Q_{(\epsilon)}]}{\epsilon}+[R, Q_{(\epsilon)}]$ with $[R, Q_{(\epsilon)}]\in O(\epsilon)$, and so 
\[\sigma_{aff}\left(\frac{\tilde{P}_{(\epsilon)}\ast_{CF} Q_{(\epsilon)}- Q_{(\epsilon)}\ast_{CF} \tilde{P}_{(\epsilon)}}{\epsilon}\right)=\sigma_{aff}\left(\frac{P_{(\epsilon)}\ast_{CF} Q_{(\epsilon)}- Q_{(\epsilon)}\ast_{CF} P_{(\epsilon)}}{\epsilon}\right).\]

As a corollary one sees that $\mathbb{S}(\g,\h,\lambda)$  is a Poisson subalgebra of $\left(S(\g)/S(\g)\h_{\lambda}\right)^{\h}$ and furthermore if $H^0_{(\epsilon)}(\h^{\bot}_{\lambda},d^{(\epsilon)}_{\h^{\bot}_{\lambda},\mathfrak{q}})$ is commutative then $\mathbb{S}(\g,\h,\lambda)$ is a Poisson commutative subalgebra of $\left(S(\g)/S(\g)\h_{\lambda}\right)^{\h}$. We recall the Vergne polarization for the sequence: Let $\g$ be a real (nilpotent) Lie algebra and let $S$ be a flag of ideals $\g_0=0\subset\g_1\subset\cdots\subset\g_n=\g$ of $\g$. Let $l\in\g^{\ast}$  and  $B_l$ be the corresponding Kostant form. Set $\g_i(l)=\{X\in\g_i /\forall Y\in\g_i,\;\; B_l(X,Y)=0\}$. Then $\mathfrak{b}_{l,S}:=\bigoplus_{i=1}^n \g_i(l)$ is a polarization called "the Vergne polarization" ( ~\cite{BER} ch. IV, Prop 1.1.2). 
\newtheorem{sistana}[ok]{Proposition}
\begin{sistana}\label{commute}
Let $\g$ be a real nilpotent Lie algebra, $\h$ and $\lambda$ as usual. If $(S(\g)/S(\g)\h_{\lambda})^{\h}$ is Poisson commutative then $H^0_{(\epsilon)}(\h^{\bot}_{\lambda},d^{(\epsilon)}_{\h^{\bot}_{\lambda},\mathfrak{q}})$ is commutative and $H^0_{(\epsilon=1)}(\h^{\bot}_{\lambda},d^{(\epsilon=1)}_{\h^{\bot}_{\lambda},\mathfrak{q}})$ is commutative.
\end{sistana}
\textit{Proof.} Suppose that $(S(\g)/S(\g)\h_{\lambda})^{\h}$ is Poisson commutative. Since $\h$ is nilpotent,  $\left(Frac(S(\g)/S(\g)\h_{\lambda})\right)^{\h}$, is also commutative (see ~\cite{BAT} ch.1, $\mathcal{x}$ 1.3.3 for this result of Dixmier). Thus by ~\cite{CG3}, the $H-$ orbits are generically lagrangian, that is $\exists\mathcal{O}\subset\lambda +\h^{\bot}$ non-empty, Zariski-open subset such that $\forall l\in\mathcal{O},\;\dim(H\cdot l)=\frac{1}{2}\dim(G\cdot l)$. In other words, $\forall l\in\mathcal{O},\;\h$ is lagrangian with respect to $l$. Also, since $\g$ is nilpotent there is a polarization $\mathfrak{b}_l$ for each $l\in\mathcal{O}$ (e.g the Vergne polarization). Thus we can apply Proposition \ref{manychars} and construct a character $\gamma_{CT}^l$ of the algebra $(U(\g)/U(\g)\h_{\lambda})^{\h}$ for each $l\in\mathcal{O}$. Suppose that $H^0_{(\epsilon)}(\h^{\bot}_{\lambda},d^{(\epsilon)}_{\h^{\bot}_{\lambda},\mathfrak{q}})$ isn't commutative. Then $\exists P_1,P_2\in H^0_{(\epsilon)}(\h^{\bot}_{\lambda},d^{(\epsilon)}_{\h^{\bot}_{\lambda},\mathfrak{q}})$, such that $Q:=P_1\ast_{CF,\epsilon}P_2-P_2\ast_{CF,\epsilon}P_1 \in H^0_{(\epsilon)}(\h^{\bot}_{\lambda},d^{(\epsilon)}_{\h^{\bot}_{\lambda},\mathfrak{q}})$ and $Q\neq 0$. Set $Q=Q_0+\epsilon Q_1+\cdots$ with $Q_0$ non homogeneous term and $Q_0\in\mathbb{R}[-f+\h^{\bot}]^{\h}$. We can also suppose without loss of generality that $Q_0\neq 0$ (if necessary we divide $Q$ by an appropriate power of $\epsilon$). Then $\gamma_{CT}^l(Q)(l)=Q\ast_{1,\mathfrak{b}} 1\in\mathbb{R}[\epsilon]$. 
The term without $\epsilon$ in $Q\ast_{1,\mathfrak{b}} 1$ is clearly $Q_0(-l)$ so $\forall l\in~\mathcal{O},\;\;\gamma_{CT}^l(Q)(l)=~0\Rightarrow \forall l\in\mathcal{O},\;\;Q_0(-l)=0\Rightarrow Q_0=0$ which is a contradiction by construction of $Q$. Thus $Q=0$ and $P_1\ast_{CF,\epsilon}P_2=P_2\ast_{CF,\epsilon}P_1$.  So  $H^0_{(\epsilon)}(\h^{\bot}_{\lambda},d^{(\epsilon)}_{\h^{\bot}_{\lambda},\mathfrak{q}})$ is commutative and of course $H^0_{(\epsilon=1)}(\h^{\bot}_{\lambda},d^{(\epsilon=1)}_{\h^{\bot}_{\lambda},\mathfrak{q}})$ is then also commutative. $\diamond$

Let $T$ be a central new variable and define $\g_T:=\g\oplus <T>$ and $\h_T:=\h\oplus <T>$ such that $\dim(\g_T)=\dim(\g)+1$. Set also $U_T(\g):=U(\g_T)$ to be the U.E.A of $\g_T$ and $U(\g_T)\h^T_{\lambda}$ to be the ideal of $U(\g_T)$ generated by  $\h_{\lambda}^T=<H+T\lambda(H),\;H\in\h>$. Following (\ref{Koornwider}), we have $(U(\g_T)/U(\g_T)\h^T_{\lambda})^{\h_T}\simeq\mathbb{D}(\g_T,\h_T,\lambda)=:\mathbb{D}_T(\g,\h,\lambda)$. One can also define as $\mathbb{D}_{(T=1)}(\g,\h,\lambda):= \mathbb{D}_T(\g,\h,\lambda)/<T-1>$ the corresponding specialization algebra. Recall that in Theorem 6.8 of \cite{BAT1}, we showed that $H^0_{(\epsilon=1)}(\h^{\bot}_{\lambda},d^{(\epsilon=1)}_{\h^{\bot}_{\lambda},\mathfrak{q}})\simeq\mathbb{D}_{(T=1)}(\g,\h,\lambda)$. The lagrangian condition (\ref{lagnow}) implies that $H^0_{(\epsilon=1)}(\h^{\bot}_{\lambda},d^{(\epsilon=1)}_{\h^{\bot}_{\lambda},\mathfrak{q}})\simeq\mathbb{D}_{(T=1)}(\g,\h,\lambda)$ is commutative as a corollary of Proposition \ref{commute}. 

\subsection{The character from representation theory.}

\textbf{Construction of a real character.} We first define a distribution on $G/B$. Let $\g$ be a nilpotent Lie algebra, $\h$ a subalgebra and $\lambda$ a real character of $\h$. Let $f\in\g^{\ast}$ s.t $f|{\h}=\lambda$ and $\mathfrak{b}$ be a polarization of $f$ with corresponding group $B$. In this section all characters are real, and if there's no confusion we'll use the previous notations for them. We denote as $C^{\infty}_c(G,B,\chi_f)$ the $C^{\infty}$ functions $\psi$ with compact support satisfying the equivariance condition for $g\in G,X\in\mathfrak{b}, \exp(X)=b$, $\psi(gb)=e^{-f(X)}\psi(g)$ and define a distribution $\alpha(f)$ by the formula for $\psi\in C^{\infty}_c(G,B,\chi_f)$, 

\begin{equation}\label{real distribution}
<\alpha(f),\psi>:=\int_{H/H\cap B}\psi(h)e^{f(H)}\mathrm{d}_{H/H\cap B}(h).
\end{equation}

Since the space $<H+f(H),\;H\in\h>$ acts by zero on $\alpha(f)$, the algebra $ (U_{\mathbb{C}}(\g)/U_{\mathbb{C}}(\g)\h_{\lambda})^{\h}$ acts on $\alpha(f)$. In \cite{BAT1}, Corollary 6.5 we proved that $D_{(T=1)}(\g,\h,\lambda)\hookrightarrow \left((U(\g)/U(\g)\h_{\lambda})^{\h}\right)$, denoted as $\mathfrak{i}_{(\epsilon=1)}$ the injective map $\mathfrak{i}_{(\epsilon=1)}:\;H^0_{(\epsilon=1)}(\h_{\lambda}^{\bot},d^{(\epsilon=1)}_{\h_{\lambda}^{\bot},\mathfrak{q}})\hookrightarrow(U(\g)/U(\g)\h_{\lambda})^{\h}$ and as $\mathcal{P}_{(t=1)}\left((U(\g)/U(\g)\h_{\lambda})^{\h}\right)$ the values at $t=1$ of polynomial in $t$ families $t\longrightarrow u_t\in \left((U(\g)/U(\g)\h_{t\lambda})^{\h}\right)$.

Let $\mathfrak{z}$ be the center of $\g$. The proof of Theorem \ref{fujisimple} was reduced to the case $\mathfrak{z}\subset \h,\;ker(f)\cap\mathfrak{z}=\{0\}$. In that case $\dim(\mathfrak{z})=1$, and if $\mathfrak{z}=<Z>$, and then there are $X,Y\in ker(f)\cap\g$ such that $[X,Y]=Z$,  and $\g=<X>\oplus\g_0$ where $\g_0:=\{W\in\g/[W,Y]=0\}$ (see ~\cite{DIX1}, $\mathcal{x}$ 4.7.7). Also, Theorem \ref{fujisimple} was proved under the hypothesis that the character $f$ was unitary. This condition was only used to establish that in the case $\h\subset\g_0$, we have $\left((U(\g)/U(\g)\h_{if})^{\h}\right)=\left((U(\g_0)/U(\g_0)\h_{if})^{\h}\right).$ The triquantization setting we used so far, concerned only the real case. Thus, before we proceed it's necessary to establish Theorem \ref{fujisimple} for a real character. 

\newtheorem{opli}[ok]{Theorem}
\begin{opli}\label{fujireal}
Let $\g$ be a finite dimensional Lie algebra, $\h\subset\g$ a subalgebra and $\lambda$ a character of $\h$. Suppose that generically  the $H-$ orbits are lagrangian in $\lambda+\h^{\bot}$. Then for regular $f\in \lambda+\h^{\bot}$ s.t $\dim(\h\cdot f)=\frac{1}{2}\dim(\g\cdot f)$ and for $A\in D_{(T=1)}(\g,\h,\lambda)$, $\mathrm{d}\tau_f(A)(\alpha(f))$ is a multiple of $\alpha(f)$, and there is defined a character $\lambda_{(T=1)}^f:\;D_{(T=1)}(\g,\h,\lambda)\longrightarrow \mathbb{R}$ such that
\begin{equation}\label{rep character}
\mathrm{d}\tau_f(A)(\alpha(f))=\lambda_{(T=1)}^f(A)\alpha(f).
\end{equation}
\end{opli}

\textit{Proof.}
The proof is based on the proof of ~\cite{FUJI2}, Theorem 1. However, we'll need some modifications in the arguments, very important for the rest of the paper. We can suppose that $\mathfrak{z}\subset\h$. Indeed, let $K$ be the corresponding group of the Lie subalgebra $\mathfrak{k}:=\h+\mathfrak{z}$. Since $\mathfrak{z}\subset\mathfrak{b}$,  for $\psi\in C^{\infty}_c(G,B,\chi_f)$ we have 
\begin{equation}\label{sec a}
<\alpha(f),\psi>=\int_{H/H\cap B}\psi(h)\chi_{\lambda}\mathrm{d}_{H/H\cap B}=\int_{HZ/(HZ)\cap B}\psi(h)\chi_{\mu}\mathrm{d}_{H/H\cap B}=\int_{K/K\cap B}\psi(k)\chi_{\mu}\mathrm{d}_{K/K\cap B},
\end{equation}
where $\mu:=f|_{(\h+\mathfrak{z})}$ and $Z$ is the corresponding group of $\mathfrak{z}$. Then for $X\in\mathfrak{k}$, $\mathrm{d}\tau_f (X) \alpha(f)= - \mu(X)\alpha(f)$ and thus  $(U_{\mathbb{C}}(\g)/U_{\mathbb{C}}(\g)\mathfrak{k}_\mu)^\mathfrak{k}$ acts on $\alpha(f)$. Also, there is a natural projection $\mathrm{J}_{\mathfrak{k}}:\; \left(U(\g)/U(\g)\h_{\lambda}\right)^{\h}\longrightarrow  \left(U(\g)/U(\g)\mathfrak{k}{\mu}\right)^{\mathfrak{k}}$ and so a natural projection $\mathrm{j}_{\mathfrak{k}}:\; D_{(T=1)}(\g,\h,\lambda)\longrightarrow D_{(T=1)}(\g,\mathfrak{k},\mu)$. Thus if $A\in D_{(T=1)}(\g,\h,\lambda)$ and $\overline{A}=\mathrm{j}_{\mathfrak{k}}(A)$ is its image in $D_{(T=1)}(\g,\mathfrak{k},\mu)$, then $\mathrm{d}\tau_f(A)(\alpha(f))=\mathrm{d}\tau_f(\overline{A})(\alpha(f))$. Since $f$ is supposed to be generic and $f\in\mu+\mathfrak{k}^{\bot}$ we can apply the induction hypothesis and there is defined a real character $\tilde{\lambda}^f_{(T=1)}:\;D_{(T=1)}(\g,\mathfrak{k},\mu)\longrightarrow \mathbb{R}$ which coincides with $\lambda^f_{(T=1)}$. So we can reduce the proof of the Theorem to the case $\mathfrak{z}\subset\h$. Furthermore, let $\mathfrak{z}^{'}:=\mathfrak{z}\cap ker(f)$. Setting $\g^{'}:=\g/\mathfrak{z}^{'}, \h^{'}:=\h/\mathfrak{z}^{'}$, we have that the conditions of the Theorem are satisfied for $\g^{'},\h^{'}, f^{'}:=f|_{\g^{'}}$. Also, $D_{(T=1)}(\g,\h,\lambda)=D_{(T=1)}(\g^{'},\h^{'},\lambda)$ and $\alpha(f)$ is $H^{'}-$ semi-invariant. The character $\lambda^{f'}_{(T=1)}:\;D_{(T=1)}(\g^{'},\h^{'},\lambda)\longrightarrow \mathbb{R}$ coincides with $\lambda^f_{(T=1)}$ so we can further reduce the proof of the Theorem to the case $\mathfrak{z}\cap ker(f)=\{0\}$. Then we distinguish two cases: First, suppose that $\h\not\subset \g_0$. Let $X\in\h\cap ker (f)$ s.t $\g=\g_0\oplus X$. Set also $\h_0=\h\cap \g_0\;\; (\dim \h_0 <\dim \h)$. Let $u_T\in U(\g_T)$. Then $u_T$ can be written as $u_T=\sum_k v_T^{(k)} X^k$ where $v_T\in U((\g_0)_T)$. Since $X\in\h\cap ker(f)$ and $U(\g_T)\h_T^{\lambda}$ is the left ideal of $U(\g_T)$ generated by elements of the form $\{H+\lambda(H)T,\;H\in\h\}$, we have for $\h^0_T:=\h_0\oplus<T>$, 
 \begin{equation}\label{uts}
 \left(U(\g_T)/U(\g_T)\h_T^{\lambda}\right)^{\h_T}\hookrightarrow  \left(U((\g_0)_T)/U((\g_0)_T)\h_T^{0\lambda}\right)^{\h_T}\subset  \left(U((\g_0)_T)/U((\g_0)_T)\h_T^{0\lambda}\right)^{\h^0_T}.
 \end{equation}
 Specializing (\ref{uts}) at $T=1$ we get for the corresponding algebras of polynomial families that 
 \begin{equation}\label{spec}
 D_{(T=1)}(\g,\h,\lambda)\hookrightarrow D_{(T=1)}(\g_0,\h_0,\lambda_0).
 \end{equation}
 
 Thus for $u_T\in \left((U(\g_T)/U(\g_T)\h_T^{\lambda})^{\h_T}/<T-1>\right)=D_{(T=1)}(\g,\h,\lambda)$, and for its image through (\ref{spec}) $u_T^0\in \left((U((\g_0)_T)/U((\g_0)_T)\h_T^{0\lambda})^{\h_T}/<T-1>\right)=D_{(T=1)}(\g_0,\h_0,\lambda_0)$ with $u_T\equiv u_T^0 mod[U(\g_T)\h_T^{\lambda}]$, we have $\lambda_{(T=1)}^f(u_T)=\lambda_{(T=1)}^{f_0}(u_T^0)$, following the corresponding computation at Theorem \ref{fujisimple}. For the proof of the case $\h\subset\g_0$ we need the following intermediate lemma:  Fix a supplementary space $\mathfrak{q}$ s.t $\g=\h\oplus\mathfrak{q}$. Let $u\in U_{\mathbb{C}}(\g)$. Let $\pi_{it\lambda}$ denote the canonical projection $\pi_{it\lambda}:\;U_{\mathbb{C}}(\g)\longrightarrow U_{\mathbb{C}}(\g)/U_{\mathbb{C}}(\g)\h_{it\lambda}$. We write $u\mapsto \pi_{it\lambda}(u)$. For $H\in\h$, let $\mathrm{d}_H:\;S_{\mathbb{C}}(\mathfrak{q})\longrightarrow S_{\mathbb{C}}(\mathfrak{q})$ denote the linear map defined by $Q\mapsto \beta^{-1}\left(\pi_{it\lambda}([H,\beta(Q)])\right)$. For $n\in \mathbb{N}$ we denote as $U_0=\mathbf{K}$ and for $n>0$ as $U_n(\g)$, the vector subspace of $U(\g)$ generated by the products $X_1X_2\cdots X_k$ with $X_1, X_2, \ldots ,X_k \in \g$ basis variables and $k\leq n$.

\newtheorem{stoi}[ok]{Lemma}
\begin{stoi}
For any $n\in \mathbb{N}$, the algebra $(U_n(\g)/U_n(\g)\h_{it\lambda})^{\h}$ depends rationally on $t$ through the symmetrization map $\beta_t$.
\end{stoi}
\textit{Proof.} Let $\{H_1,\ldots, H_t\}$ be a basis of $\h$ s.t $\lambda(H_t)=1$ and $\{H_1,\ldots,H_{t-1}\}$ is a basis of $ker \lambda \cap \h$. Let $\{Q_1,\ldots,Q_r\}$ a basis of $\mathfrak{q}$. To prove this Lemma, we check the kernels of the maps $\mathrm{d}_{H_i}$, that is the elements $Q^{\alpha}=Q_1^{\alpha_1}\cdots Q_r^{\alpha_r}\in S_{\mathbb{C}}(\mathfrak{q})$ s.t $\pi_{it\lambda}([H_i,\beta(Q^{\alpha})])=0$. We have that $\forall i$, $[H_i,Q^{\alpha}]=\sum_{\gamma,\delta}c^i_{\gamma\delta}Q^{\delta}H^{\gamma}$. For $z\in\mathbb{C}$, we have $mod U_{\mathbb{C}}(\g)\h_{z\lambda}$ that 
\[[H_i,Q^{\alpha}]\equiv\sum_{\delta}c^i_{\delta}Q^{\delta}+\sum_{\gamma=(0,\ldots,0,\gamma_t),\delta}c^i_{\gamma\delta}Q^{\delta}(-z)^{\gamma_t}.\]

Applying $\beta^{-1}$ on $\pi_{z\lambda}([H_i,\beta(Q^{\alpha})])$ in order to argue in $S_{\mathbb{C}}(\mathfrak{q})$, we see that the linear applications $\mathrm{d}_{H_i}$ correspond to matrices with polynomial coefficients, and thus we can do calculations in the field of rational fractions $\mathbb{C}(z)$. Our algebra of invariants corresponds to the common kernel of $\mathrm{d}_{H_i}$ which depends in a rational way on $z$, with $z$ generic. Indeed, let's first write $\mathrm{d}_{H_i}(z)$ meaning the dependance on $z$. Examining the dimension of the kernels $ker(\mathrm{d}_{H_i}(z))$ for generic $z$, (actually proving that $\forall i,j,\;\dim_{\mathbb{C}(z)}(ker \mathrm{d}_{H_i}\cap ker \mathrm{d}_{H_j})=\dim_{\mathbb{C}}(ker \mathrm{d}_{H_i}(z)\cap ker\mathrm{d}_{H_j}(z))$ and that $[ker\mathrm{d}_{H_i}\cap ker\mathrm{d}_{H_j}](z)=ker \mathrm{d}_{H_i}(z)\cap ker \mathrm{d}_{H_j}(z)$), we conclude that the linear maps $\mathrm{d}_{H_i}$ depend rationally on generic values of $z$ to prove this intermediate Lemma. $\diamond$

Before we continue recall that the orbit $H\cdot f$ is said to be saturated with respect to $\g_0^{\bot}$ iff $H\cdot f+ \g_0^{\bot}=H\cdot f$. We now continue the proof of Theorem \ref{fujireal} for the case $\h\subset\g_0$. In this case, the condition of Corwin-Greenleaf (see (1) of equations (2.7) in ~\cite{FLMM}) holds for the character $it\lambda$ and we have 
\begin{equation}\label{ratiot}
\left(U_{\mathbb{C}}(\g)/U_{\mathbb{C}}(\g)\h_{it\lambda}\right)^{\h}=\left(U_{\mathbb{C}}(\g_0)/U_{\mathbb{C}}(\g_0)\h_{it\lambda}\right)^{\h}.
\end{equation}
This equation depends rationally on $it$, $t\in\mathbb{R}^{\ast}$. So if (\ref{ratiot}) holds for $it$, $t\in\mathbb{R}^{\ast}$, it holds also for $t$ in a Zariski-open subset of $\mathbb{R}$ and we write 
\begin{equation}\label{ratiot1}
\left(U_{\mathbb{C}}(\g)/U_{\mathbb{C}}(\g)\h_{t\lambda}\right)^{\h}=\left(U_{\mathbb{C}}(\g_0)/U_{\mathbb{C}}(\g_0)\h_{t\lambda}\right)^{\h}, 
\end{equation}
but we cannot conclude that a similar equation holds for the algebra $\left(U_{\mathbb{C}}(\g)/U_{\mathbb{C}}(\g)\h_{\lambda}\right)^{\h}$. This is the difference with respect to the proof of Theorem \ref{fujisimple} which was about a unitary character. To overcome this setback, we can use polynomial families $t\mapsto u_t\in  \left(U_{\mathbb{C}}(\g_0)/U_{\mathbb{C}}(\g_0)\h_{t\lambda}\right)^{\h}$. This will allow us to continue the argument and it explains at the same time why the character at the Theorem's statement is defined for $D_{(T=1)}(\g,\h,\lambda)$ contrary to Theorem \ref{fujisimple} where the result holds for $\left(U_{\mathbb{C}}(\g)/U_{\mathbb{C}}(\g)\h_{\lambda}\right)^{\h}$. We argue as follows: Again for $T$ such that $[T,\g]=0$, set $\g_T:=\g\oplus\mathbb{R}<T>$,  $(\g_0)_T:=\g_0\oplus\mathbb{R}<T>$, and $\h_T=\{H+T\lambda(H),\;H\in\h\}\subset \g_T$. Since here $\h\subset \g_0$, we also have $\h_T\subset (\g_0)_T$. Finally set $U(\g_T)$ to be the U.E.A of $\g_T$ and $U((\g_0)_T)$ the U.E.A of $(\g_0)_T$. Recall also that if $\h\subset\g_0$, then $\forall f\subset \mathcal{O}\subset \lambda+\h^{\bot}$, the orbits $H\cdot f$ are saturated with respect to $\g_0^{\bot}$. Using this fact we'll prove that $(U((\g_0)_T)/U((\g_0)_T)\h_{T})^{\h_T}= (U(\g_T)/U(\g_T)\h_{T})^{\h_T}$. First we have that $l\in\h_T^{\bot}\Leftrightarrow l= f+\mu T^{\ast}$ for $f\in\h_{\mu\lambda}^{\bot}$ (recall that abusively we wrote $\h_{\lambda}^{\bot}=-\lambda+\h^{\bot}$). Indeed for $H\in\h$ we have that $l\in\h_T^{\bot}\Leftrightarrow l(H+T\lambda(H))=0\Leftrightarrow l(H)+l(T)\lambda(H)=0\Leftrightarrow l=f+\mu T^{\ast},\;f\in\h_{\mu \lambda}^{\bot}$. Thus for $l\in\h_T^{\bot},\;X\in\g_T$, we also have that $X\in\h_T(l)\Leftrightarrow\;\forall H'\in\h_T,\;l([X,H']=0\Leftrightarrow \;\forall H\in\h,\;f([X,H])=0\Leftrightarrow X\in\h(f)\oplus<T>$. Since for the data $\g,\h,\lambda$ and $f\in\mathcal{O}\subset\h_{\lambda}^{\bot}$, the $H\cdot f$ orbits are saturated with respect to $\g_0^{\bot}$, we have $\dim(\h(f))=\dim(\h(f_0))-1$, for  $f\in\mathcal{O}$ (see (2.7)-(2.8) of \cite{FLMM}). Set $\mathcal{O}_T:=\cup_{\mu} \mu \mathcal{O}$. The set $\mathcal{O}_T$ is Zariski open since $\mathcal{O}$ is itself Zariski open.  Thus we have that for $l\in\mathcal{O}_T$,  $\dim(\h_T(l))=\dim(\h_T(l|_{(\g_0)_T}))-1$ and so the $H_T-$ orbits are saturated with respect to $(\g_0)_T^{\bot}\subset \g_T^{\ast}$. Theorem 5.2 from ~\cite{FLMM} along with the equivalences (2.7), (2.8) from ~\cite{FLMM}, essentially prove that for generic $f\in\mathcal{O}\subset \h^{\bot}_{\lambda}$,
\begin{equation}\label{Theoremagain}
\mathcal{D}(\g_0,\h,\lambda_0)=\mathcal{D}(\g,\h,\lambda)\Leftrightarrow H\cdot f \;\;\textlatin{ are saturated with respect to}\;\;\g_0^{\bot}.
\end{equation}
Applying Theorem (\ref{Theoremagain}), for $\g_T,\h_T$ we get (again for $U_{\mathbb{C}}(\g)$)

\begin{equation}\label{forwi}
(U_{\mathbb{C}}((\g_0)_T)/U_{\mathbb{C}}((\g_0)_T)\h_{T})^{\h_T}=(U_{\mathbb{C}}(\g_T)/U_{\mathbb{C}}(\g_T)\h_{T})^{\h_T}.
\end{equation}
 Specializing (\ref{forwi}) at $T=1$, we have that the algebras $\left((U_{\mathbb{C}}((\g_0)_T)/U_{\mathbb{C}}((\g_0)_T)\h_{T})^{\h_T}/<T-1>\right)$ and $\left((U_{\mathbb{C}}(\g_T)/U_{\mathbb{C}}(\g_T)\h_{T})^{\h_T}/<T-1>\right)$, are equal and correspond to the same polynomial family of operators, namely
\begin{equation}\label{same ds}
 D_{(T=1)}(\g_0,\h,\lambda_0)=D_{(T=1)}(\g,\h,\lambda).
 \end{equation}
Recall that $D_{(T=1)}(\g,\h,\lambda)\hookrightarrow \left(U_{\mathbb{C}}(\g)/U_{\mathbb{C}}(\g)\h_{\lambda}\right)^{\h}$ so in this case we can provide the Theorem only on the level of $D_{(T=1)}(\g,\h,\lambda)$. The character in this case is calculated accordingly with Theorem \ref{fujisimple}. $\diamond$
 
 \subsection{The Theorem on characters}
\newtheorem{cla}[ok]{Theorem}
\begin{cla}\label{central}
Let $\g$ be a finite dimensional nilpotent Lie algebra, $\h\subset \g$ a subalgebra, $\lambda$ a character of~$\h$. Let also $P\in  H^0_{(\epsilon=1)}(\h_{\lambda}^{\bot},d^{(\epsilon=1)}_{\h_{\lambda}^{\bot},\mathfrak{q}})$ and $u\in D_{(T=1)}(\g,\h,\lambda)$ such that $u= \mathfrak{i}_{(\epsilon=1)}(P)$. Then for generic $f\in\lambda+\h^{\bot}$, there is a pair $(\mathfrak{b}_f, \mathfrak{q}_f)$ such that 
\begin{equation}\label{equal characters}
T_1\circ\overline{\beta}^{-1}_{\mathfrak{q}_f,(\epsilon)}(P)|_{\epsilon=1}(-f)=\lambda^f_{(T=1)}(u).
\end{equation}
\end{cla}
Before its proof we have to examine the conditions of Theorems \ref{fujisimple} and \ref{fujireal} and check if the special case for which those Theorems were proved, can be applied also in biquantization without loss of generality. There is an important difference in the initial conditions needed to construct a character with biquantization techniques and harmonic analysis. In the first case, the construction can be given "pointwise". That is for a chosen $f\in\lambda+\h^{\bot}$ such that the lagrangian condition holds, we create a character of $H^0_{(\epsilon=1)}(\h^{\bot}_f,d^{(\epsilon=1)}_{\h^{\bot}_f,\mathfrak{q}})$ like in Proposition \ref{ct}. If moreover the algebra $(U(\g)/U(\g)\h_{\lambda})^{\h}$ is commutative then we can extend the construction and create such a character for $H^0_{(\epsilon=1)}(\h^{\bot}_{\lambda},d^{(\epsilon=1)}_{\h^{\bot}_{\lambda},\mathfrak{q}})$ (recall that $H^0_{(\epsilon=1)}(\h^{\bot}_{\lambda},d^{(\epsilon=1)}_{\h^{\bot}_{\lambda},\mathfrak{q}})$ is smaller than $\left(U(\g)/U(\g)\h_{\lambda}\right)^{\h}$) at a generic $f\in\lambda+\h^{\bot}$. In the case of harmonic analysis, the statement of ~\cite{FUJI2} Theorem 1 (which is Theorem \ref{fujisimple} here) supposes the finite multiplicity condition for the representation~$\tau_{\lambda}$. Thus the character of $(U(\g)/U(\g)\h_{\lambda})^{\h}$ is constructed there for generic $f$ satisfying the lagrangian condition. But Theorem \ref{central} refers to a (maybe) smaller algebra $D_{(T=1)}(\g,\h,\lambda)\subset \left(U(\g)/U(\g)\h_{\lambda}\right)^{\h}$. Its proof uses double induction on $\dim(\g)$ and $\dim(\h)$. For this we proceed in three steps: First we prove that one can suppose that $\mathfrak{z}\subset\h$. This is done by showing that if instead of $\h$ we use the Lie subalgebra $\mathfrak{k}:=\h+\mathfrak{z}$ from the biquantization side we always calculate the same character. Once this is established, we prove that one can also suppose that $\mathfrak{z}\cap ker(f)=\{0\}$. Again we are concentrated in the biquantization setting. The last part is induction on $\dim(\g)$ where we distinguish the cases $\h\subset\g_0$ and $\h\not\subset \g_0$.

\textbf{First step: The hypothesis $\mathfrak{z}\subset \h$.} Let $\mathfrak{k}=\h+\mathfrak{z}$ and $K\subset G$ be its corresponding Lie subgroup. Note that the following Lemma is a pointwise result.
\newtheorem{strei}[ok]{Lemma}
\begin{strei}\label{firstlevel}
Let $\g$ be a nilpotent Lie algebra, $\mathfrak{z}$ the center of $\g$, $\h\subset\g$ a subalgebra, $\lambda$ a character of $\h$, and $f\in\lambda+\h^{\bot}$. Let $\mathfrak{b}$ be a polarization with respect to $f$. Let also $\g=\h\oplus\mathfrak{q}_{\mathfrak{b}}$ be a decomposition of $\g$ where $\mathfrak{q}_{\mathfrak{b}}$ satisfies the transversal condition for the polarization $\mathfrak{b}$. Let $\mathfrak{z}=\mathfrak{z}\cap\h\oplus\mathfrak{z}\cap\mathfrak{q}_{\mathfrak{b}}$ and set $\mathfrak{z}_{\mathfrak{q}_{\mathfrak{b}}}:=\mathfrak{z}\cap\mathfrak{q}_{\mathfrak{b}}$. Finally let $\mathfrak{k}=\h\oplus\mathfrak{z}_{\mathfrak{q}_{\mathfrak{b}}}$ and $\mathfrak{q}_{\mathfrak{b}}=\mathfrak{z}_{\mathfrak{q}_{\mathfrak{b}}}\oplus V$, i.e $V$ is a supplementary of $\mathfrak{z}_{\mathfrak{q}_{\mathfrak{b}}}$ in $\mathfrak{q}_{\mathfrak{b}}$.  Then 
\[P\in H^0_{(\epsilon)}(\h^{\bot}_{\lambda},d^{(\epsilon)}_{\h^{\bot}_{\lambda},\mathfrak{q}_{\mathfrak{b}}})\Rightarrow P|_{-f+\mathfrak{k}^{\bot}}\in H^0_{(\epsilon)}(\mathfrak{k}^{\bot}_f,d^{(\epsilon)}_{\mathfrak{k}^{\bot}_f,V}).\]
\end{strei}

\textit{Proof.} Let $P=\sum c^{\alpha} z^{\alpha}P_{\alpha},\;z^{\alpha}\in S(\mathfrak{z}_{\mathfrak{q}_{\mathfrak{b}}})[\epsilon],\;P_{\alpha}\in S(V)[\epsilon]$. Set $\overline{P}:=P|_{-f+\mathfrak{k}^{\bot}}$. In the equations defining the reduction algebra (\cite{BAT1} after Lemma 3.4) 
we saw that only two kinds of Kontsevich graphs arise, namely $\mathcal{B}-$ type and $\mathcal{BW}-$ type graphs. Let's formulate the reduction equations for $ H^0_{(\epsilon)}(\mathfrak{t}^{\bot}_f,d^{(\epsilon)}_{\mathfrak{t}^{\bot}_f,V})$. For this, let $D_n$ be the differential operator corresponding to a possible $\mathcal{B}_n-$ type graph. Then $D_n(P)=\sum c^{\alpha} z^{\alpha}D_n(P_{\alpha})$. Examining the possible coloring in such a $\mathcal{B}-$ type graph, the edges deriving $\overline{P}$ carry variables from $V$. Furthermore, no vertex except the root of the $\mathcal{B}-$ type graph can belong to $\mathfrak{z}_{\mathfrak{q}_{\mathfrak{b}}}$. Let $\overline{D_n},\overline{P}$ be the respective objects defined on the new supplementary $V$, in the decomposition $\g=\mathfrak{k}\oplus V$. Thus checking all cases of coloring for the decomposition $\g=\h\oplus\mathfrak{q}_{\mathfrak{b}}$ we see that then

\begin{equation}\label{reduction center}
\overline{D_n(P)}=\overline{D_n}(\overline{P}).
\end{equation}

With the same reasoning the Lemma holds for $D_n$ being an operator coming from a $\mathcal{BW}_n-$ type graph. Finally, since for the differential $d^{(\epsilon)}_{\h^{\bot}_{\lambda},\mathfrak{q}_{\mathfrak{b}}}$ we have $d^{(\epsilon)}_{\h^{\bot}_{\lambda},\mathfrak{q}_{\mathfrak{b}}}=\sum_nD_n$,  we have by (\ref{reduction center})
\[P\in H^0_{(\epsilon)}(\h^{\bot}_{\lambda},d^{(\epsilon)}_{\h^{\bot}_{\lambda},\mathfrak{q}_{\mathfrak{b}}})\Rightarrow \sum_i\epsilon^id^{(i)}_{\h^{\bot}_{\lambda},\mathfrak{q}_{\mathfrak{b}}}(P)=0\Rightarrow \sum_i\epsilon^id^{(i)}_{\mathfrak{k}^{\bot}_l,V}(\overline{P})=0\Rightarrow \overline{P}\in  H^0_{(\epsilon)}(\mathfrak{k}^{\bot}_f,d^{(\epsilon)}_{\mathfrak{k}^{\bot}_f,V}).\;\;\diamond\]

\newtheorem{srpo}[ok]{Lemma}
\begin{srpo}
Let $\g,\h,\lambda,\mathfrak{z}$ be as before, $V$ be a vector subspace such that $\g=(\h+\mathfrak{z})\oplus V$ and $\mathfrak{q}$ such that $\g=\h\oplus\mathfrak{q}$. Let $f\in \lambda+\h^{\bot}$ and suppose the lagrangian condition holds for $f$. Let $\gamma_{CT}^{\mathfrak{z}}:  (U_{(\epsilon)}(\g)/U_{(\epsilon)}(\g)(\h+\mathfrak{z})_f)^{\h+\mathfrak{z}}\longrightarrow~\mathbb{R}[\epsilon],\; u\mapsto \overline{T}_1\circ\overline{\beta}^{-1}_{V,(\epsilon)}(u)$ and $\gamma_{CT}:  (U_{(\epsilon)}(\g)/U_{(\epsilon)}(\g)\h_f)^{\h}\longrightarrow \mathbb{R}[\epsilon],\; u\mapsto \overline{T}_1\circ\overline{\beta}^{-1}_{\mathfrak{q},(\epsilon)}(u)$ be characters constructed as in Theorem \ref{ctb}. Then 
\begin{equation}\label{equal first case}
\overline{T}_1\circ\overline{\beta}_{V,(\epsilon)}^{-1}(\overline{u})=\overline{T}_1\circ\overline{\beta}_{\mathfrak{q},(\epsilon)}^{-1}(u)|_{-f+(\h+\mathfrak{z})^{\bot}}.
\end{equation}
\end{srpo}

\textit{Proof.}
 Let $P=\sum_{\alpha} z_{\alpha}Q_{\alpha} \in S(\mathfrak{q})$ with $z_{\alpha}\in S(\mathfrak{z}_{\mathfrak{q}})[\epsilon], Q_{\alpha}\in S(V)[\epsilon]$.
Let $u=\overline{\beta}_{\mathfrak{q},_{(\epsilon)}}(P)+U_{(\epsilon)}(\g)\h_f$. We have $\overline{\beta}_{\mathfrak{q},(\epsilon)}(P)=\sum_{\alpha} z_{\alpha}\beta_{(\epsilon)}(Q_{\alpha})$. Set $\overline{u}=\sum z_{\alpha}\beta(Q_{\alpha}).$ Then from Lemma \ref{firstlevel}, the two characters in (\ref{equal first case}) are well defined and equal. $\diamond$

We then consider the characters of the specialized algebras $\gamma_{CT}^{(\epsilon=1)}:\;H^0_{(\epsilon=1)}(\h^{\bot}_{\lambda},d^{(\epsilon=1)}_{\h^{\bot}_{\lambda},\mathfrak{q}_{\mathfrak{b}}})\longrightarrow \mathbb{R}$ and $\gamma_{CT}^{\mathfrak{z},(\epsilon=1)}:\;H^0_{(\epsilon=1)}(\mathfrak{k}^{\bot}_f,d^{(\epsilon=1)}_{\mathfrak{k}^{\bot}_f,V})\longrightarrow \mathbb{R}$ for which we'll have $\gamma_{CT}^{(\epsilon=1)}(u)=\gamma_{CT}^{\mathfrak{z},(\epsilon=1)}(\overline{u})$ (recall that the operators in Theorem \ref{ctb} are composed of $\mathcal{W}-$ type graphs so intervention of $\mathfrak{z}$ doesn't alter these computations). On the harmonic analysis side, we define a character of $\mathfrak{k}$ as $\mu:=f|_{(\h+\mathfrak{z})}$. There is a morphism $\delta:\;\left(U(\g_T)/U(\g_T)\h_T^{\lambda}\right)^{\h_T}\longrightarrow \left(U(\g_T)/U(\g_T)\mathfrak{k}_T^{\mu}\right)^{\mathfrak{k}_T}$ and by specialization at $T=1$ we get a morphism 
\[\delta_{(T=1)}:\; D_{(T=1)}(\g,\h,\lambda)\longrightarrow D_{(T=1)}(\g,\mathfrak{k},\mu)\subset \left(U(\g)/U(\g)\mathfrak{k}_{\mu}\right)^{\mathfrak{k}}.\]

Then by Theorems \ref{fujisimple} and \ref{fujireal} we get that for $u\in D_{(T=1)}(\g,\h,\lambda)$, $\overline{u}$ its image in $D_{(T=1)}(\g,\mathfrak{k},\mu)$ and $\overline{\lambda}_{(T=1)}^f$ the character of Theorem \ref{fujireal} for $D_{(T=1)}(\g,\mathfrak{k},\mu)$, it is $\overline{\lambda}^f_{(T=1)}(\overline{u})=\lambda_{(T=1)}^f(u)$. So the constructions and computations from both sides match. Thus if Theorem \ref{central} is true for $\g$, the subalgebra $\mathfrak{k}=\h+\mathfrak{z}\subset \g$, and a polarization $\mathfrak{b}$, then one can take 
\[\mathfrak{b}':=\mathfrak{b}\cap(\h+\mathfrak{z}_{\mathfrak{q}_{\mathfrak{b}}})\oplus\mathfrak{b}\cap V=(\mathfrak{b}\cap \h+\mathfrak{z}_{\mathfrak{q}_{\mathfrak{b}}})\oplus \mathfrak{b}\cap V=\mathfrak{b}\cap\h\oplus\mathfrak{b}\cap(\mathfrak{z}_{\mathfrak{q}_{\mathfrak{b}}}+ V),\]
and prove the same fact for the subalgebra $\h$ and the polarization $\mathfrak{b}'$ simply by taking for supplementary the space $\mathfrak{q}=\mathfrak{z}_{\mathfrak{q}_{\mathfrak{b}}} + V$. We thus saw that the character from biquantization can be calculated supposing $\mathfrak{z}\subset\h$ without loss of generality. For the analytic setting this was known from Theorems \ref{fujisimple} and \ref{fujireal}.

\textbf{Second step: The hypothesis $\mathfrak{z}\cap ker(\lambda)=\{0\}$.}
\newtheorem{stpl3}[ok]{Lemma}
\begin{stpl3}
Let $\g$ be a nilpotent Lie algebra, $\h$ a subalgebra, $\lambda$ a character of $\h$ and $\mathfrak{z}$ the center of $\g$. Suppose that $\mathfrak{z}\subset\h$. Let $\mathfrak{q}$ be such that $\g=\h\oplus\mathfrak{q}$, and set $\mathfrak{z} ':=\mathfrak{z}\cap ker(\lambda)$. Let $\g':=\g/\mathfrak{z}'$, $\h':=\h/\mathfrak{z}'$ and let $\mathfrak{q}'$ be such that $\g'=\h'\oplus\mathfrak{q}'$. Then
\[H^0_{(\epsilon)}(\h^{\bot}_{\lambda},d^{(\epsilon)}_{\h^{\bot}_{\lambda},\mathfrak{q}})\simeq H^0_{(\epsilon)}((\h')^{\bot}_{\lambda},d^{(\epsilon)}_{(\h')^{\bot}_{\lambda},\mathfrak{q}'}).\]
\end{stpl3}
\textit{Proof.} Since $\mathfrak{z}'\subset \h$ and $\mathfrak{q}'$ is such that $\g'=\h'\oplus\mathfrak{q}'$, we have $(\mathfrak{q}')^{\ast}\simeq \mathfrak{q}^{\ast}$. Recall that $\mathfrak{q}^{\ast}\simeq \h^{\bot}$.
Furthermore $(\h')^{\bot}:=\{X\in \g/ \forall Y\in\h',\;\lambda([X,Y])=0\}$ and since again $\mathfrak{z}'\subset\mathfrak{z}$ and $\h=\h'+\mathfrak{z}'$, $(\h')^{\bot}=\h^{\bot}$.
Thus in the following  diagram all the arrows are vector space isomorphisms.
\begin{center}
\begin{tabular}{ccc}
$(\h')^{\bot}$ & $\rightarrow$ & $\h^{\bot}$\\
$\downarrow$ &   & $\downarrow$\\
$(\mathfrak{q}')^{\ast}$ & $\rightarrow$ & $\mathfrak{q}^{\ast}$\\
\end{tabular} 
\end{center}
Check the differentials $d^{(\epsilon)}_{\h^{\bot}_{\lambda},\mathfrak{q}}$ and $d^{(\epsilon)}_{(\h')^{\bot}_{\lambda},\mathfrak{q}'}$. Since $\h'=\h/\mathfrak{z}'$, the graphs in $d^{(\epsilon)}_{(\h')^{\bot}_{\lambda},\mathfrak{q}'}$ are the same with the graphs in the differential $d^{(\epsilon)}_{\h^{\bot}_{\lambda},\mathfrak{q}'}$. Then by ~\cite{CT} a change of supplementary space for $\h^{\bot}_{\lambda}$ gives isomorphic reduction algebras so $H^0_{(\epsilon)}(\h^{\bot}_{\lambda},d^{(\epsilon)}_{\h^{\bot}_{\lambda},\mathfrak{q}'})\simeq H^0_{(\epsilon)}(\h^{\bot}_{\lambda},d^{(\epsilon)}_{\h^{\bot}_{\lambda},\mathfrak{q}})$ and so $H^0_{(\epsilon)}(\h^{\bot}_{\lambda},d^{(\epsilon)}_{\h^{\bot}_{\lambda},\mathfrak{q}})\simeq H^0_{(\epsilon)}((\h')^{\bot}_{\lambda},d^{(\epsilon)}_{(\h')^{\bot}_{\lambda},\mathfrak{q}'})$. $\diamond$

Thus, under the condition $\mathfrak{z}'=\mathfrak{z}\cap\ker(f)=\{0\}$ there is no loss of generality for the character calculated with biquantization. For the harmonic analysis setting this is also true from Theorems \ref{fujisimple} and \ref{fujireal}. Thus we can now proceed to the last step in the proof of Theorem \ref{central} under the conditions $\mathfrak{z}\subset \h$, $\mathfrak{z}\cap\ker(f)=\{0\}$.

\textbf{Third step: Induction on $\dim\g$.} 

Under these conditions we have $\dim(\mathfrak{z})=1$. Let $\mathfrak{z}=\mathbb{R}Z$. Then there are $X,Y\in \g\cap ker(f)$ such that $[X,Y]=Z,\; f(Z)=1$ (~\cite{DIX1} $\mathcal{x}$ 4.7.7). If we also set $\g_0:=\{W\in\g/[W,Y]=0\}$, then $\g=\g_0 + <X>$ ($\g_0$ is an hyperplane and codimension 1 ideal). In his proof, Fujiwara supposes that $\mathfrak{b}\subset\g_0$. As we'll see this implies a particular choise of $\mathfrak{q}_{\mathfrak{b}}$. We proceed by induction. For $\dim(\g)=1$ the Theorem is true since $\h=\g$ and $(U(\g)/U(\g)\h_{\lambda})^{\h}=\mathbb{R}$. Suppose by induction that the Theorem is true for all nilpotent Lie algebras $\g$ such that the lagrangian condition holds and $\dim(\g)=n-1$. In the case of $\g_0$, the induction hypothesis holds for $\g_0,\h,\lambda_0=\lambda|_{\g_0}$. For the final step of the induction, we consider two cases:

\textbf{Case} $\mathbf{\h\subset\g_0}$.  
Since $\h\subset \g_0$, let $V$ be such that $\g_0=\h\oplus V$. Applying the induction hypothesis to $\g_0,\h,V$, one can take $\mathfrak{q}:=V\oplus <X>$ as a new supplementary in order to prove the claim for $\g,\h,\mathfrak{q}$. In biquantization terms we have that 
\begin{equation}\label{two special}
H^0_{(\epsilon=1)}(\h_{\lambda}^{\bot},d^{(\epsilon=1)}_{\h_{\lambda}^{\bot},V})=H^0_{(\epsilon=1)}(\h_{\lambda}^{\bot},d^{(\epsilon=1)}_{\h_{\lambda}^{\bot},\mathfrak{q}}),
\end{equation}
where the first is considered as a subalgebra of $S(\g_0)$ and the second as a subalgebra of $S(\g)$. Indeed, (\ref{two special}) holds: \cite{BAT1}, Theorem 6.8 in view of (\ref{forwi}) says that

\begin{equation}\label{oneto}
H^0_{(\epsilon=1)}(\h_{\lambda}^{\bot},d^{(\epsilon=1)}_{\h_{\lambda}^{\bot},V})=\mathbb{D}_{(T=1)}(\g_0,\h,\lambda_0)=\mathbb{D}_{(T=1)}(\g,\h,\lambda)=H^0_{(\epsilon=1)}(\h_{\lambda}^{\bot},d^{(\epsilon=1)}_{\h_{\lambda}^{\bot},\mathfrak{q}}).
\end{equation}
Finally since $\mathfrak{i}_{(\epsilon=1)}(H^0_{(\epsilon=1)}(\h_{\lambda}^{\bot},d^{(\epsilon=1)}_{\h_{\lambda}^{\bot},V}))=\mathcal{P}_{(t=1)}\left((U(\g_0)/U(\g_0)\h_{\lambda})^{\h}\right)$ and $ \mathfrak{i}_{(\epsilon=1)}(H^0_{(\epsilon=1)}(\h_{\lambda}^{\bot},d^{(\epsilon=1)}_{\h_{\lambda}^{\bot},\mathfrak{q}}))=\mathcal{P}_{(t=1)}\left((U(\g)/U(\g)\h_{\lambda})^{\h}\right)$, 
we can write $\mathfrak{i}_{(\epsilon=1)}(H^0_{(\epsilon=1)}(\h_{\lambda}^{\bot},d^{(\epsilon=1)}_{\h_{\lambda}^{\bot},\mathfrak{q}}))=\mathcal{P}_{(t=1)}\left((U(\g)/U(\g)\h_{\lambda})^{\h}\right)=\newline
=\mathcal{P}_{(t=1)}\left((U(\g_0)/U(\g_0)\h_{\lambda})^{\h}\right)= \mathfrak{i}_{(\epsilon=1)}(H^0_{(\epsilon=1)}(\h_{\lambda}^{\bot},d^{(\epsilon=1)}_{\h_{\lambda}^{\bot},V})$. Thus we proved that all spaces and algebras behave in the same way in our induction on the harmonic analysis and the biquantization side.

We prove now that the calculation of characters also coincides. Indeed, since $\mathfrak{b}\subset\g_0$, the distribution $\alpha(f)$ can now be considered as a distribution on $G_0/B$ and so by (\ref{same ds}), for the characters $\lambda_{(T=1)}^f:\;D_{(T=1)}(\g,\h,\lambda)\longrightarrow\mathbb{R}$, and $\lambda_{(T=1)}^{f'}:\;D_{(T=1)}(\g_0,\h,\lambda)\longrightarrow \mathbb{R}$ we have $\lambda_{(T=1)}^f=\lambda_{(T=1)}^{f'}$. Indeed, in the biquantization side, equality (\ref{two special}) suggests that every $P\in H^0_{(\epsilon=1)}(\h_{\lambda}^{\bot},d^{(\epsilon=1)}_{\h_{\lambda}^{\bot},\mathfrak{q}})$ can be written as an element of $S(V)$. Thus since $\h\subset\g_0,\;\mathfrak{b}\subset\g_0$, at the corner of the biquantization diagram of $-f+\mathfrak{b}^{\bot}$ and $-f+\h^{\bot}$ where the character $\gamma_{CT}:\;H^0_{(\epsilon)}(\h_{\lambda}^{\bot},d^{(\epsilon)}_{\h_{\lambda}^{\bot},\mathfrak{q}})\longrightarrow \mathbb{R}[\epsilon],\;P\mapsto T_2^L|_{-f+(\h+\mathfrak{b})^{\bot}}(P)$ is calculated, we see that all components involved in the computation are in $\g_0$. Thus if $\gamma^{(\epsilon=1)'}_{CT}:\;H^0_{(\epsilon=1)}(\h_{\lambda_0}^{\bot},d^{(\epsilon=1)}_{\h_{\lambda_0}^{\bot},V})\longrightarrow \mathbb{R}$ is the character computed for $\g_0,\h,\lambda_0$, we have that for $P\in  H^0_{(\epsilon=1)}(\h_{\lambda}^{\bot},d^{(\epsilon=1)}_{\h_{\lambda}^{\bot},\mathfrak{q}})$,  $\gamma^{(\epsilon=1)'}_{CT}(P)=\gamma^{(\epsilon=1)}_{CT}(P)$. Apply the induction hypothesis for $\g_0,\h,\lambda_0$ to prove the Theorem in the case $\h\subset\g_0$.

  \textbf{Case} $\mathbf{\h\not\subset \g_0}$. Set again $\h_0:=\g_0\cap \h$. Suppose that for $\dim(\g_0)=n-1$, $\h_0\not\subset \h$, the claim holds for $\g_0,\h_0,\lambda_0:=\lambda|_{\g_0}$. Let $U(\g_0)\h_{0\lambda}$ denote the ideal of $U(\g_0)$ generated by the family of elements $<H+\lambda_0(H),\;H\in\h_0>$.
Since $\h\not\subset \g_0$, we can use the same supplementary to pass from $\g_0$ to $\g$. That is we write $\g=\h\oplus\mathfrak{q}$ and $\g_0=\h_0\oplus\mathfrak{q}$. Regarding the reduction algebras $H^0_{(\epsilon)}(\h_{\lambda}^{\bot},d^{(\epsilon)}_{\h_{\lambda}^{\bot},\mathfrak{q}})\;,\; H^0_{(\epsilon)}(\h_{0\lambda}^{\bot},d^{(\epsilon)}_{\h_{0\lambda}^{\bot},\mathfrak{q}})$ we have now that $H^0_{(\epsilon)}(\h_{\lambda}^{\bot},d^{(\epsilon)}_{\h_{\lambda}^{\bot},\mathfrak{q}})\subset H^0_{(\epsilon)}(\h_{0\lambda}^{\bot},d^{(\epsilon)}_{\h_{0\lambda}^{\bot},\mathfrak{q}})$.  Indeed the condition $\h\not\subset\g_0$ implies that there is $Y\in\h$ s.t $Y\not\subset\g_0$ and so a decomposition $\g=\mathfrak{q}\oplus\h_0\oplus<Y>$. Then since $\g_0=\mathfrak{q}\oplus\h_0$ is an ideal of $\g$, we have $[\g_0,\g]\subset\g_0$ which means that $[Y,\mathfrak{q}]\subset\g_0$ and so $d^{(\epsilon)}_{\h_{\lambda}^{\bot},\mathfrak{q}}$ contains all the possible graphs in $d^{(\epsilon)}_{\h_{0\lambda}^{\bot},\mathfrak{q}}$ and additionally those that the variable $Y$  is associated to the edge $e_{\infty}$. Thus $H^0_{(\epsilon)}(\h_{\lambda}^{\bot},d^{(\epsilon)}_{\h_{\lambda}^{\bot},\mathfrak{q}})\subset H^0_{(\epsilon)}(\h_{0\lambda}^{\bot},d^{(\epsilon)}_{\h_{0\lambda}^{\bot},\mathfrak{q}})$, and specializing at $\epsilon=1$ we have $H^0_{(\epsilon=1)}(\h_{\lambda}^{\bot},d^{(\epsilon=1)}_{\h_{\lambda}^{\bot},\mathfrak{q}})\subset H^0_{(\epsilon=1)}(\h_{0\lambda}^{\bot},d^{(\epsilon=1)}_{\h_{0\lambda}^{\bot},\mathfrak{q}})$. With respect to the computation of the character from harmonic analysis in Theorem 9, we have in this case that for $u_T\in \left((U(\g_T)/U(\g_T)\h_T^{\lambda})^{\h_T}/<T-1>\right)=D_{(T=1)}(\g,\h,\lambda)$, and for its image $u_T^0\in \left((U((\g_0)_T)/U((\g_0)_T)\h_T^{0\lambda})^{\h_T}/<T-1>\right)=D_{(T=1)}(\g_0,\h_0,\lambda_0)$ with $u_T\equiv u_T^0 mod[U(\g_T)\h_T^{\lambda}]$, $\lambda_{(T=1)}^f(u_T)=\lambda_{(T=1)}^{f_0}(u_T^0)$. Respectively in the biquantization setting, the supplementary $\mathfrak{q}$ as we said is now in $\g_0$. Since also $\mathfrak{b}\subset\g_0$, the calculations for  $\gamma_{CT}:\;H^0_{(\epsilon)}(\h_{\lambda}^{\bot},d^{(\epsilon)}_{\h_{\lambda}^{\bot},\mathfrak{q}})\longrightarrow \mathbb{R}[\epsilon]$
are taking place in $\g_0$. So if $\gamma_{CT}'':\;H^0_{(\epsilon)}(\h_{0\lambda}^{\bot},d^{(\epsilon)}_{\h_{0\lambda}^{\bot},\mathfrak{q}_0})\longrightarrow \mathbb{R}[\epsilon]$, is the character computed for $\g_0,\h_0,\lambda_0$, let $\gamma^{(\epsilon=1)''}_{CT}:\;H^0_{(\epsilon=1)}(\h_{0\lambda}^{\bot},d^{(\epsilon=1)}_{\h_{0\lambda}^{\bot},\mathfrak{q}_0})\longrightarrow \mathbb{R}$ and $\gamma^{(\epsilon=1)}_{CT}:\;H^0_{(\epsilon=1)}(\h_{\lambda}^{\bot},d^{(\epsilon=1)}_{\h_{\lambda}^{\bot},\mathfrak{q}})\longrightarrow \mathbb{R}$ be the corresponding characters defined on the specializations. Then we get that for $P\in  H^0_{(\epsilon=1)}(\h_{\lambda}^{\bot},d^{(\epsilon=1)}_{\h_{\lambda}^{\bot},\mathfrak{q}})\subset  H^0_{(\epsilon=1)}(\h_{0\lambda}^{\bot},d^{(\epsilon)}_{\h_{0\lambda}^{\bot},\mathfrak{q}})$, we have $\gamma^{(\epsilon=1)''}_{CT}(P)=\gamma^{(\epsilon=1)}_{CT}(P)$. $\diamond$

It is important here to note that the construction of the character $\gamma_{CT}(P)=T_1\circ\overline{\beta}_{\mathfrak{q}_f}^{-1}(P)$ depends on the choice of the supplementary space $\mathfrak{q}_f$. Since this space is chosen with respect to the polarization $\mathfrak{b}_f$ the character also depends on $\mathfrak{b}_f$. This is in contrast with the representation theory character $u\mapsto\lambda_f(u)$ which, as was noted before, is independent of the polarization used. So to prove the equality (\ref{equal characters}), we have to choose $\mathfrak{b}_f,\mathfrak{q}_f$. 

\section{On an example of H.Fujiwara.}
In this section we discuss an example of H.Fujiwara to show the importance of the transversality condition between the polarization and the supplementary of $\h$, and provide a character formula. For the choise of $\mathfrak{q}$ made originally by Fujiwara, we show that there is no polarization in transversal position 
\begin{equation}\label{transversality condition}
\mathfrak{b}=\mathfrak{b}\cap\mathfrak{q} + \mathfrak{b}\cap \h
\end{equation}
with respect to this specific $\mathfrak{q}$. Then we verify that the $\ast-$ product for $S(\mathfrak{q})^{\h}$ in this case is the Moyal product and make evident that the symmetrization map $\beta_{\mathfrak{q}}:\;S(\mathfrak{q})^{\h}\longrightarrow(U_{\mathbb{C}}(\g)/U_{\mathbb{C}}(\g)\h_{\lambda})^{\h}$ is not an isomorphism of algebras for this choise of $\mathfrak{q}$. Afterwards we compute a transversal choice $\mathfrak{q}_V$ of $\mathfrak{q}$ for the Vergne polarization. The symmetrization operator $\beta_{\mathfrak{q}_V}$ is subject to change if we change $\lambda$ and as a consequence $\beta_{\mathfrak{q}_V}$ isn't an isomorphism of algebras. Instead we suggest that $\beta_{\mathfrak{q}_V}$ is an isomorphism of algebras only when it is evaluated on $l\in\lambda +\h^{\bot}$. This is verified on the generators of $S(\mathfrak{q})^{\h}$. The operators $T_1,T_2$ are trivial in this example so $\beta_{\mathfrak{q}_V}^{-1}$ gives a character for $(U_{\mathbb{C}}(\g)/U_{\mathbb{C}}(\g)\h_{\lambda})^{\h}$. All computations that are left out can be found at ~\cite{BAT} $\mathcal{x}$ 5.5.
\subsection{Description of the example.}
Consider the real, nilpotent Lie algebra $\g$, generated by elements ${X,U,V,E,Z}$, with non-zero brackets 
\[ [U,V]=E, [X,U]=V, [X,V]=Z.\]
Let $\lambda=E^{\ast}\in\g^{\ast}$, and consider the Lie subalgebra $\h=\mathbb{R}X\oplus \mathbb{R}E$ of $\g$.
Let $\mathfrak{z}(U_{\mathbb{C}}(\g))$ be the center of the envelopping algebra $U_{\mathbb{C}}(\g)$, and take the elements 
\[A=2UZ-V^2 , \;W=A-2EX\]
of $U_{\mathbb{C}}(\g)$. $A$ is $\h-$ invariant and also $[\g,W]=0$ with $W \in \mathfrak{z}(U_{\mathbb{C}}(\g))$. For $l\in \lambda + \h^{\bot}\simeq \mathbb{R}^3$, we have $\g(l) = \mathfrak{z}(\g) + \mathbb{R}(X - l(Z)U + l(V)V)$. As $EX \in U_{\mathbb{C}}(\g)\h_{\lambda} = <H+  \lambda(H),\;H\in\h>$, $A$ and $W$ are in the same coset $modU_{\mathbb{C}}(\g)\h$ so represent the same element of $\mathbb{D}(\g,\h,\lambda)\simeq (U_{\mathbb{C}}(\g)/U_{\mathbb{C}}(\g)\h_{\lambda})^{\h}$. The algebra $\mathbb{D}(\g,\h,\lambda)$ is commutative as the $H\cdot l$ orbits are lagrangian. Indeed, both $\h,\;\h^l:=\{A\in\g/\forall B\in \h,\;\;l([A,B])=0\}$ are isotropic subspaces with respect to $l$. Let now $\mathfrak{q}=<Z,U,V>$. For this choise of supplementary space, we have a symmetric space $\g=\h\oplus\mathfrak{q}$. To calculate the algebra $S(\mathfrak{q})^{\h}$ we first calculate the coadjoint orbits of $H$, the associated to $\h$ connected, nilpotent Lie group. With the appropriate parametrizations and eliminations one eventually finds that  
\begin{equation}\label{invariants}
S(\mathfrak{q})^{\h}\simeq\mathbb{C}[Z,2ZU-V^2].
\end{equation}
For the choise $\mathfrak{q}=<Z,U,V>$, we have by direct computation $A^2= 4U^2Z^2 - 4Z\beta(UV^2)+ V^4= \beta(A^2).$

\subsection{First choise of $\mathfrak{q}$.}
Let now $\bar{A}=2UZ-V^2$ denote the element of $S_{\mathbb{C}}(\mathfrak{q})^{\h}$. Since the choise $\mathfrak{q}=<Z,U,V>$ makes $(\h,\mathfrak{q})$ a symmetric pair, we know that $H^0_{(\epsilon)}(\h^{\bot}_{\lambda},d^{(\epsilon)}_{\h^{\bot}_{\lambda},\mathfrak{q}})\simeq (S_{\mathbb{C}}(\mathfrak{q})^{\h}[\epsilon],\ast_{CF})$. In ~\cite{CT}, ch.3, it is proved that in the case of symmetric spaces, the $\ast_{CF}-$product coincides with the product given by the function $E(X,Y)$, where  \[E_{\lambda}(X,Y):=\exp({\lambda(H(X,Y))})E(X,Y),\;\;B_{\pi}(\exp(X),\exp(Y))=E_{\lambda}(X,Y)\exp({X+Y}),\]
$\lambda\in\h^{\ast},\;X,Y\in\mathfrak{q},$
and $B_{\pi}$ is the product defined in ~\cite{CT}, $\mathcal{x}$1.4.1. Furthermore, in this example, we may see $\mathfrak{q}$ as a symplectic space, suggesting that actually the $\ast-$ product in $S(\mathfrak{q})^{\h}$ is nothing else than the Moyal product $\ast_M$. Calculating this product and powers of $A$ we take $U\ast_M V^2=UV^2 +\{U,V^2\}=UV^2+M_{12}\frac{\partial U}{\partial U}\cdot \frac{\partial V^2}{\partial V}=UV^2 -V$ and $V^2\ast_M U= V^2U +\{V^2,U\}=V^2U+ M_{21}\frac{\partial V^2}{\partial V}\cdot \frac{\partial U}{\partial U}=V^2U+V$.

Thus, in $S(\mathfrak{q})^{\h}$ we compute:
\begin{equation}\label{barA^2 in moyal}
\bar{A}\ast_M \bar{A}= (2ZU-V^2)\ast_M (2ZU-V^2)=4Z^2U^2-2Z(U\ast_M V^2)-2Z( V^2\ast_M U) +V^4=\bar{A}^2
\end{equation}

In order to prove an isomorphism of algebras $(S_{\mathbb{C}}(\mathfrak{q})^{\h},\ast_M)\simeq (S_{\mathbb{C}}(\mathfrak{q})^{\h},\cdot)$, problems arise with powers $\bar{A}^k$, and $A^k,$ for $k\geq 3$.
To attack the problem from both sides, we want to compute $\bar{A}^3$ for the Moyal product. For this, we use the graph description of $\ast_M$. Putting $\bar{A}^2$, $\bar{A}$ on the horizontal axis, we run over the possible nonzero graphs
to construct $\ast_M$. After this check and taking in mind the coefficients coming from the Lie bracket and the formula of $\ast_M$, one finds
\begin{equation}\label{equality products}
\bar{A}\ast_M\bar{A}\ast_M\bar{A}=\bar{A}^2\ast_M \bar{A}= \bar{A}^2\cdot\bar{A}+ \frac{\partial^2}{\partial U^2}\bar{A}^2\cdot\frac{\partial^2}{\partial V^2}\bar{A}=\bar{A}^3- \frac{1}{8}\cdot 8Z^2\cdot 2=\bar{A}^3-2Z^2,
\end{equation}
where in the left hand-side we have the Moyal product in $S_{\mathbb{C}}(\mathfrak{q})^{\h}$, and in the right-hand side the ordinary product in $S_{\mathbb{C}}(\mathfrak{q})^{\h}$.
Thus, it is obvious that $\beta_{\mathfrak{q}}:\;(S_{\mathbb{C}}(\mathfrak{q})^{\h},\cdot)\longrightarrow (U_{\mathbb{C}}(\g)/U_{\mathbb{C}}(\g)\h_{\lambda})^{\h},$ for $\mathfrak{q}=<Z,U,V>$ is not an isomorphism of algebras. As an immediate consequence, for $l\in \lambda +\h^{\bot}$, the map 
\[(U_{\mathbb{C}}(\g)/U_{\mathbb{C}}(\g)\h_{\lambda})^{\h}\ni A\mapsto \beta_{\mathfrak{q}}^{-1}(A)(l)\]
isn't a morphism of algebras, i.e a character of $(U_{\mathbb{C}}(\g)/U_{\mathbb{C}}(\g)\h_{\lambda})^{\h}$.

We claim that to construct such a character, one needs to select $\mathfrak{q}$ in a transversal way with respect to a given polarization $\mathfrak{b}$ and the subalgebra $\h$ and then evaluate on $l$. In this example the particular choise of $\mathfrak{q}$ made so far is not transversal to any polarization. 

Indeed, let $\mathfrak{q}=\mathbb{R}<U,V,Z>$. This is an $\h-$ invariant supplementary space of $\h$, isomorphic to $\g/\h$. As $\dim(\g)=5, \dim(\g(l))=3$, a polarization $\mathfrak{b}$ has to be of dimension $\frac{1}{2}(\dim(\g)+ \dim(\g(l)) = 4$. By definition of the polarization it has to be $\g(l)\subseteq \mathfrak{b}$ and by checking the dimensions, $\g(l)$ will be a codimension 1 subalgebra of $\mathfrak{b}$, so an ideal of $\mathfrak{b}$. So finding a polarisation $\mathfrak{b}$, is the same thing as to find a vector $Y \in \g$ such that $Y\notin \g(l)$ and $[Y,\g(l)]\subseteq \g(l)$.

The condition $l([Y,\g(l)])=0$ will be always satisfied by the definition of $\g(l)$. It turns out that we have only one possibility for $Y$. Indeed, since $\forall l\in \lambda+\h^{\bot},\;\;\g(l)=<E,Z>\oplus\mathbb{R}(X-l(Z)U+l(V)V)$, is of codimension 1 in $\mathfrak{b}$, it is also an ideal of $\mathfrak{b}$. Thus, it should be $[Y,\g(l)]\subset\g(l)$ and in particular, $[Y,K]\subset\g(l)$, where $K=X-l(Z)U$. Assume that $Y=aU+bV$. Then $[Y,K]=-aV-bZ-l(Z)bE$ which has to be in $\g(l)$. This means $a=0$ and so $Y=V$.
Thus, the condition (\ref{transversality condition}) is not satisfied and the choise $\mathfrak{q}=\mathbb{R}<Z,U,V>$ for the supplementary of $\h$ is not transversal to any polarization $\mathfrak{b}$. This explains the fact that $\beta_{\mathfrak{q}}$ is not an isomorphism of algebras as we already calculated in (\ref{equality products}).
 
\subsection{A transversal choise of $\mathfrak{q}$.}
Let $l\in \lambda+ \h^{\bot}$. For such an $l$ and choosing the flag $S$ as \[\g_1 =\mathbb{R}E, \;\g_2 =\mathbb{R}E \oplus \mathbb{R} Z,\;\g_3=\mathbb{R}E \oplus \mathbb{R} Z \oplus \mathbb{R}V,\;\g_4=\mathbb{R}E \oplus \mathbb{R} Z \oplus \mathbb{R}V \oplus \mathbb{R}X,\;\g_5=\mathbb{R}E \oplus \mathbb{R} Z \oplus \mathbb{R}V \oplus \mathbb{R}X \oplus \mathbb{R}U,\]
then one can compute the Vergne polarization with respect to these data and get
 \begin{equation}\label{q depends on l}
\mathfrak{b}_V(l,S)=<E,Z,V,K>,\;\; \mathfrak{q}_V=<Z,V,K>.
\end{equation}
for $l \in \mathcal{O}=\{f\in\lambda+h^{\bot}/f(Z)\neq 0\}$,  and $K=X-l(Z)U+l(V)V$. This is a transversal choise of $\mathfrak{q}$ and we also compute that $T_{1,\mathfrak{q}_V}=1$. An important remark here is that $\mathfrak{b}_V$ is the only possible polarization in this example. 

\subsection{Comparison of the two choises.}
We denote as $\mathfrak{q}_l$, instead of $\mathfrak{q}_V$, the transversal choise suggested by the Vergne polarisation. This makes sense, as $\mathfrak{q}_l$ in (\ref{q depends on l}) depends on $l\in\lambda+\h^{\bot}$.
Let also $\bar{\beta}_{\mathfrak{q}}:\;S_{\mathbb{C}}(\mathfrak{q})\longrightarrow (U_{\mathbb{C}}(\g)/U_{\mathbb{C}}(\g)\h_{\lambda})$ and $\bar{\beta}_{\mathfrak{q}_l}:\;S_{\mathbb{C}}(\mathfrak{q}_l)\longrightarrow (U_{\mathbb{C}}(\g)/U_{\mathbb{C}}(\g)\h_{\lambda})$ denote the corresponding symmetrization maps for each choise of $\mathfrak{q}$ \footnote{We can't take directly the invariant subspaces here as $\mathfrak{q}_l$ isn't $\h-$ invariant.}.
Let  $\tau_l$ be the irreducible representation associated to the orbit $G\cdot l$ by the Kirillov map and let $\alpha_l$ be the element of $(\mathcal{H}_l^{-\infty})^{H,\chi_{\lambda}}$ defined in (\ref{distribution}). As  $W \in \mathfrak{z}U_{\mathbb{C}}(\g)$, we have that $\beta^{-1}(W^3) = \beta^{-1}(W)^3,$ and we calculate the action of  $A$ and $W$ on $\alpha_l$ through the induced by $\g$ representation of $U_{\mathbb{C}}(\g)$ on $(\mathcal{H}_l^{-\infty})^{H,\chi_{\lambda}}$:
\[ \mathrm{d}\tau^{-\infty}_l(A)\alpha_l=\mathrm{d}\tau^{-\infty}_l(W)\alpha_l=\beta^{-1}(A)(il)\alpha_l=(l(V)^2-2l(Z)l(U))\alpha_l.\]
Also,
\[\mathrm{d}\tau^{-\infty}_l(A^3)\alpha_l=\mathrm{d}\tau^{-\infty}_l(W^3)\alpha_l=\beta^{-1}(W^3)(il)\alpha_l=(l(V)^2-2l(Z)l(U))^3\alpha_l.\]
However, if we calclate directly in $S(\g)$, we get 
\begin{equation}\label{extra term}
\beta^{-1}(A^3)=(2UZ-V^2)^3-2Z^2E^2,
\end{equation}

that is 
\[\mathrm{d}\tau^{-\infty}_l(A^3)\alpha_l\neq \beta^{-1}(A^3)(il)\alpha_l.\]

We will show now that $\bar{\beta}_{\mathfrak{q}_l}:\;S_{\mathbb{C}}(\mathfrak{q}_l)\longrightarrow (U_{\mathbb{C}}(\g)/U_{\mathbb{C}}(\g)\h_{\lambda})$ induces a character $A\mapsto \beta_{\mathfrak{q}_l}^{-1}(A)(l)$ for the algebra of invariant elements when evaluated on $l\in\lambda+\h^{\bot}$. Take again $K=X-l(Z)U$, $\mathfrak{q}_l=<Z,V,K>$, $\bar{A}=2UZ-V^2$. Then 
\[A=2UZ-V^2=2Z(\frac{X-K}{l(Z)})-V^2\;\; \textlatin{and}\;\; A\equiv [-\frac{2ZK}{l(Z)}-V^2]\quad mod\;U_{\mathbb{C}}(\g)\h_{\lambda}.\]
 Thus, $\bar{\beta}_{\mathfrak{q}_l}^{-1}(A)=-\frac{2ZK}{l(Z)}-V^2$ which is well defined as we supposed that $l(Z)\neq 0$. When evaluated on $l$, we have $\bar{\beta}_{\mathfrak{q}_l}^{-1}(A)(l)=-l(V)^2-2l(K)=-l(V)^2-2l(X)+2l(Z)l(U)=\bar{A}(l)$, and the same for $\bar{A}^2$: $\bar{\beta}_{\mathfrak{q}_l}^{-1}(A^2)(l)=\bar{A}^2(l),\;\bar{\beta}_{\mathfrak{q}_l}^{-1}(A^2)=\bar{A}^2$. We'll see now that $\bar{\beta}_{\mathfrak{q}_l}^{-1}(A^3)(l)=\bar{A}^3(l)$, while $\bar{\beta}_{\mathfrak{q}_l}^{-1}(A^3)\neq\bar{A}^3$. Indeed, this is because one has to change $\mathfrak{q}_l$ when $l(Z)$ changes.

Let $\hat{E}=E-\frac{Z}{l(Z)},\;\;\hat{K}=U-\frac{X}{l(Z)}$. Then $[\hat{K},V]=E-\frac{Z}{l(Z)}=\hat{E},\;\;[X,\hat{K}]=V=[X,U]$. Thus the Lie algebra $\g$ of the example, is isomorphic to $\hat{\g}=<\hat{E},X,\hat{K},V,Z>$ and in particular, setting $\hat{\h}=<\hat{E},X>$, $\hat{\mathfrak{q}}=\hat{\g}/\hat{\h}$, $\hat{\g}$ is a symmetric space. Then $l\in\hat{\h}^{\bot}$ since $l([\hat{K},V])=l(\hat{E})=0$ and we fall back to the case of a symmetric space without a character.

In that case, we know that the symmetrization map $\hat{\beta}_{\hat{\mathfrak{q}}}:\;S(\hat{\mathfrak{q}}_{\mathbb{C}})^{\hat{\h}}\longrightarrow (U_{\mathbb{C}}(\hat{\g})/U_{\mathbb{C}}(\hat{\g})\hat{\h})^{\hat{\h}}$ is an isomorphism of algebras. Here, $2\hat{K}Z-V^2$ is again $\hat{\h}-$ invariant and by the theory of nilpotent symmetric spaces (see ~\cite{B1}), we get $\hat{\beta}((2\hat{K}Z-V^2)^3)=(2\hat{K}Z-V^2)^3mod U_{\mathbb{C}}(\hat{\g})\hat{\h}$. 

On the other hand,
as polynomials, $2\hat{K}Z-V^2\equiv \bar{A}$ on $\h^{\bot}_{\lambda}$. Also, since $[A,X]=[A,Z]=0$ we may write 
\[\left(A-\frac{2XZ}{l(Z)}\right)^3\equiv A^3modU_{\mathbb{C}}(\g)\h_{\lambda},\;\;\left(A-\frac{2XZ}{l(Z)}\right)^3\equiv A^3modU_{\mathbb{C}}(\g)\hat{\h}.\]

By the fact that $\g\simeq\hat{\g}$ and (\ref{extra term}), we can rewrite for $\hat{\g}$ and the symmetrization map $\hat{\beta}:\;S(\hat{\g}_{\mathbb{C}})\longrightarrow~U(\hat{\g}_{\mathbb{C}})$,
\[\hat{\beta}^{-1}((2\hat{K}Z-V^2)^3)=(2\hat{K}Z-V^2)^3-2Z^2\hat{E}^2.\]
Conversely, in $U(\hat{\g}_{\mathbb{C}})$ 
\[\hat{\beta}((2\hat{K}Z-V^2)^3)-2Z^2\hat{E}^2=(2\hat{K}Z-V^2)^3=A^3mod U_{\mathbb{C}}(\hat{\g})\hat{\h}_{\lambda},\]
since $2\hat{K}Z-V^2=A-\frac{2ZX}{l(Z)}$. Finally, $modU_{\mathbb{C}}(\g)\h_{\lambda}$, we have $\hat{E}=1-\frac{Z}{l(Z)}$ and so

\begin{equation}\label{extra term final}
A^3=\beta_{\mathfrak{q}_l}(\bar{A}^3)-2Z^2(1-\frac{Z}{l(Z)})^2=\beta_{\mathfrak{q}_l}(\bar{A}^3-2Z^2(1-\frac{Z}{l(Z)})^2).
\end{equation}
One can see clearly here, that $\bar{\beta}_{\mathfrak{q}_l}^{-1}(A^3)\neq \bar{A}^3$ but $\bar{\beta}_{\mathfrak{q}_l}^{-1}(A^3)(l)= \bar{A}^3(l)$ as the second term at the right hand side of (\ref{extra term final}) disappears when we evaluate on $l$.

\subsection{Isomorphism of reduction spaces.}
We have calculated that $H^0_{(\epsilon)}(\h^{\bot}_{\lambda},d^{(\epsilon)}_{\h^{\bot}_{\lambda},\mathfrak{q}})=S(\mathfrak{q})^{\h}[\epsilon]=<Z,2UZ-V^2>[\epsilon]$. For the transversal choise $\mathfrak{b}_V,\mathfrak{q}_V$, we can simply take $K=X-l(Z)U$ as $V$ is already in the supplementary space. The question here is to compute the isomorphism
\[B_{\mathfrak{q},\mathfrak{q}_l}:\;H^0_{(\epsilon)}(\h^{\bot}_{\lambda},d^{(\epsilon)}_{\h^{\bot}_{\lambda},\mathfrak{q}})\longrightarrow H^0_{(\epsilon)}(\h_{\lambda}^{\bot},d^{(\epsilon)}_{\h_{\lambda}^{\bot},\mathfrak{q}_l}).\]

We have to describe first the $\ast -$ product in $H^0_{(\epsilon)}(\h_{\lambda}^{\bot},d^{(\epsilon)}_{\h_{\lambda}^{\bot},\mathfrak{q}_V})$. The pair $\h\oplus\mathfrak{q}_l$ is almost symmetric, since the only non-convenient bracket is the $[K,V]=Z-l(Z)E\in\mathfrak{z}(\g)$. This means that the $\ast-$ product we are looking for, will be the sum of the Moyal product and some extra terms coming from this bracket. 

We follow the approach in \cite{CT}. Finding the vector field $\upsilon$, we proceed to the solution of the differential equation
\begin{equation}\label{isomorphism of supplementaries}
\partial \mu_t=[DU_{\hat{\pi_t}}(\hat{\upsilon}),\mu_t],
\end{equation}
where $U$ is the $L_{\infty}-$ quasi-isomorphism of Kontsevich, $\pi_t=e^{t\mathrm{ad}\upsilon}\cdot \pi=\pi+t[\upsilon,\pi]+\frac{t^2}{2}[\upsilon,[\upsilon,\pi]],\;\mu_t=U(e^{\pi_t})$ and the hat denotes the partial Fourier transform. To do so, one has to check the possible graphs with non-zero contributions for the operator $DU$:

The change of basis
\[\{X,E,Z,V,U\}\longrightarrow \{X,E,Z,V,U-X/l(Z)\},\]
is non-singular for $l(Z)\neq 0$, i.e for $l\in \mathcal{O}$. Thus after the construction in ~\cite{CT} we write $\upsilon=\frac{X}{l(Z)}\partial_U$. According to the brackets from the data of the example, the Poisson vector field $\pi$ here is described by three diagrams as shown in Fig. 3.
\begin{figure}[h!]
\begin{center}
\includegraphics[width=8cm]{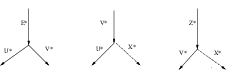}
\caption\footnotesize{The nonzero components of $\pi$}
\end{center}
\end{figure}
The computation of $[\pi,\upsilon]$ consists in gluing the graphs of the two vector fields $\upsilon,\pi$ and searching for non-zero graphs. This gives, that $[\pi,\upsilon]=\pi\circ\upsilon + \upsilon\circ\pi$ corresponds to the graph in Fig. 4.
\begin{figure}[h!]\label{braketofpi}
\begin{center}
\includegraphics[width=8cm]{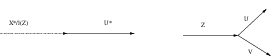}
\caption\footnotesize{At the left, the vector field $\upsilon$ and at the right the bracket $[\pi,\upsilon]$.}
\end{center}
\end{figure}

Thus, $\pi_t=\pi+t[\pi,\upsilon]$ and according with what is just said, $\pi_t$ can be expressed in graph terms. Next, one checks the graphs in the expression of the operator $DU_{\pi_t}(\upsilon)$. In general, the graph corresponding to this operator, will receive a colored $(\in \h^{\ast})$ edge from the horizontal axis, which goes either to a vertex of the graph or leaves to infinity. In the first case, there still has to be a colored edge leaving to infinity, which narrows down the open possibilities for graphs. 
\begin{figure}[h!]
\begin{center}
\includegraphics[width=9cm]{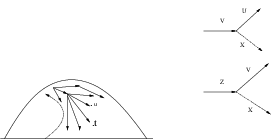}
\caption\footnotesize{At the left, how the operator $DU_{\pi_t}$ looks like and at the right two a priori nonzero graphs.}
\end{center}
\end{figure}
In our example, and after the graphs appearing in $\pi$, the ground edge leaves to infinity, so we search for contributions among the different types of graphs listed in $\mathcal{x}\;5.5.1$ of ~\cite{CT} and which we recall in Fig. 6.

\begin{figure}[h!]
\begin{center}
\includegraphics[width=8cm]{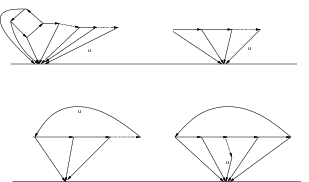}
\caption\footnotesize{Possible graphs for $DU_{\pi_t}$.}
\end{center}
\end{figure}

Using among others the argument that $U$ isn't in the derived algebra $[\g,\g]$ of $\g$, we are left with two nonzero graphs for $DU_{\pi_t}$, and after calculating their weights (equal to $\frac{1}{12}$), 
\begin{equation}\label{operator DU in example}
DU_{\pi_t}(\upsilon)=\frac{1}{12}\left(\frac{1}{l(Z)}\partial_U^3-\frac{tZ}{l(Z)^2}\partial_U^3\right).
\end{equation}
As a special but simple case, the graph in Fig. 7
\begin{figure}[h!]
\begin{center}
\includegraphics[width=7cm]{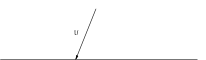}
\caption\footnotesize{Another zero contribution.}
\end{center}
\end{figure}
is also zero because when we'll later evaluate at $l\in\{\xi\;/\; \xi|_{\h}=\lambda|_{\h}\}$ we'll get $l(X)=0$ from $\upsilon=-\frac{X}{l(Z)}\partial_U,\;\lambda=E^{\ast}$. Graphically, the operator $DU_{\pi_t}$ is given by the sum at Fig. 8.

\begin{figure}[h!]
\begin{center}
\includegraphics[width=8cm]{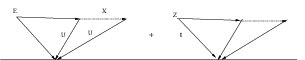}
\caption\footnotesize{The solution graphically.}
\end{center}
\end{figure}
$DU_{\pi_t}$ is a commuting family of operators and so the solution of equation (\ref{isomorphism of supplementaries}) is
\begin{equation}\label{solution of equation}
e^{\int_tDU_{\pi_t} dt}.
\end{equation}

As far as it concerns relations between $(H^0_{(\epsilon)}(\h_{\lambda}^{\bot},d^{(\epsilon)}_{\h_{\lambda}^{\bot},\mathfrak{q}}),\ast_M)\simeq (H^0_{(\epsilon)}(\h_{\lambda}^{\bot},d^{(\epsilon)}_{\h_{\lambda}^{\bot},\mathfrak{q}_l}),\ast_{M_l})$ and $(U_{\mathbb{C}}(\g)/U_{\mathbb{C}}(\g)\h_{\lambda})^{\h}$ we can write
\begin{equation}\label{relations}
\beta^{-1}((2UZ-V^2)^3-2Z^2)=\bar{A}\ast_M\bar{A}\ast_M\bar{A}=\bar{A}^3-2Z^2, \;\;
\end{equation}
\begin{equation}
 \beta_{\mathfrak{q}_l}^{-1}((2UZ-V^2)^3-2Z^2)=(2UZ-V^2)^3-2Z^2+\frac{4Z^3}{l(Z)}-\frac{2Z^4}{l(Z)^2}.
\end{equation}

Indeed, for $t=1$, \[B_{\mathfrak{q},\mathfrak{q}_l}((2UZ-V^2)^3-2Z^2)=e^{[\frac{1}{12l(Z)}(1-\frac{Z}{2l(Z)})\partial_U^3]}((2UZ-V^2)^3-2Z^2)=\]
\[((2UZ-~V^2)^3-~2Z^2)+~\frac{4Z^3}{l(Z)}-~\frac{2Z^4}{l(Z)^2}.\] Thus a first formula that one can derive is
\begin{equation}\label{two betas}
\beta_{\mathfrak{q}_l}^{-1}(A)=e^{[\frac{1}{12l(Z)}(1-\frac{Z}{2l(Z)})\partial_U^3]}\beta_{\mathfrak{q}}^{-1}(A).
\end{equation}

So after evaluation at $l$, we get 
\begin{equation}\label{re}
 \beta_{\mathfrak{q}_l}^{-1}(A^3)(l)=\bar{A}^3(l),
 \end{equation}
 and
 $\beta_{\mathfrak{q}_l}^{-1}(A)(l)=\beta_{\mathfrak{q}}^{-1}(A)(l)$ by (\ref{two betas}).
If $\phi(A)$ is the polynomial associated to $A$, this in turn implies that the map
\[(U_{\mathbb{C}}(\g)/U_{\mathbb{C}}(\g)\h_{\lambda})^{\h}\ni A\mapsto \left[l\mapsto\phi(A)(l)=\left(e^{[\frac{1}{12l(Z)}(1-\frac{Z}{2l(Z)})\partial_U^3]}\beta_{\mathfrak{q}}^{-1}(A)\right)(l)\right]\]
is polynomial.
Let $\bar{B}=P(\bar{A},\bar{Z}), B=P(Z,A)$ be two polynomial expressions in $S_{\mathbb{C}}(\mathfrak{q})^{\h}$ and $(U_{\mathbb{C}}(\g)/U_{\mathbb{C}}(\g)\h_{\lambda})^{\h}$ respectively. We will omit the bar for $Z$ in $\bar{B}$ as no confusion arises. Since we have verified the isomorphism on the generators $A,Z$ \footnote{We didn't actually verify it for Z but this is trivial as it is a central element.} one can write
\begin{equation}\label{two betas2}
\beta_{\mathfrak{q}}^{-1}(B)=e^{-\frac{1}{24Z}\partial_U^3}\bar{B},
\end{equation}
where $\beta_{\mathfrak{q}_l}^{-1}(B)=\bar{B}=P(\bar{A},Z)$. The numerator of the exponential here is justified by the fact that we haven't yet evaluated at $l$ so by (\ref{two betas}), $\frac{1}{12Z}(1-\frac{Z}{2Z})=\frac{1}{24Z}$. An easy calculation with equation (\ref{two betas2}) verifies that $\beta_{\mathfrak{q}}^{-1}(A^3)=\bar{A}^3-2Z^2$ as calculated by means of the Moyal product in (\ref{equality products}). So the final formula we get is 
\[\beta_{\mathfrak{q}}^{-1}(B)(l)=e^{-\frac{1}{24l(Z)}\partial_U^3}\bar{B}(l).\]

As a conclusion, we see that under (\ref{transversality condition}), the map $\beta_{\mathfrak{q}_l}$ becomes an isomorphism when evaluated on $l\in \lambda+\h^{\bot}$ constructing a character of $(U_{\mathbb{C}}(\g)/U_{\mathbb{C}}(\g)\h_{\lambda})^{\h}$. 

This way we verify and generalize the result of ~\cite{BL}. There the authors prove the isomorphism $(U(\g)/U(\g)\h)^{\h}\simeq (S(\g)/S(\g)\h)^{\h}$ in the case where exists a common polarization for all elements $l$. In that case, there is an obvious choise of supplementary $\mathfrak{q}_l$ suggested by the condition (\ref{transversality condition}). We generalized this result, by showing that whenever there isn't such a convenient polarization, in order to construct a character one has to choose $\mathfrak{q}_l$ transversally with respect to $\mathfrak{b}_l,\h$ and then evaluate $\beta_{\mathfrak{q}_l}$ on $l$.

On the other hand, these methods permit us to compute explicitly the pseudodifferential operator appearing in the character formula and composed with $\beta_{\mathfrak{q}}^{-1}$. This way we can write down the formula of the character, a task almost impossible otherwise.

\newtheorem{oxw}[ok]{Theorem}
\begin{oxw}
Let $\g$ be the real nilpotent Lie algebra, generated by the elements ${X,U,V,E,Z}$, with non-trivial brackets $[U,V]=E, [X,U]=V, [X,V]=Z.$  Let $\lambda=E^{\ast}\in\g^{\ast}$, and consider the Lie subalgebra $\h=\mathbb{R}X\oplus \mathbb{R}E$ of $\g$. For $\mathfrak{q}=<Z,V,U>$ and $\mathfrak{q}_l=<Z,V,X-l(Z)U+l(V)V>$ with $l(Z)\neq 0$, we have that for $u\in (U(\g)/U(\g)\h_{\lambda})^{\h}$,
\begin{equation}\label{formule}
\beta_{\mathfrak{q}_l}^{-1}(u)=e^{[\frac{1}{12l(Z)}(1-\frac{Z}{2l(Z)})\partial_U^3]}\beta_{\mathfrak{q}}^{-1}(u),
\end{equation}
and if $v$ is a polynomial in $(U(\g)/U(\g)\h_{\lambda})^{\h}$ then $\beta_{\mathfrak{q}}^{-1}(v)(l)=e^{-\frac{1}{24l(Z)}\partial_U^3}\beta^{-1}_{\mathfrak{q}_l}(v)(l)$.

The map 
\[\gamma_{CT}:\;v\mapsto \left(e^{[\frac{1}{12l(Z)}(1-\frac{Z}{2l(Z)})\partial_U^3]}\beta_{\mathfrak{q}}^{-1}(v)\right)(l)\]
is a character of $(U(\g)/U(\g)\h_{\lambda})^{\h}$.
\end{oxw}
\textbf{Acknowledgements.} This paper is the sequence of \cite{BAT1} and also comes out from the PhD thesis \cite{BAT}. The author would like to gratefully thank Charles Torossian for his inspiring guidance and supervision during these years. He would also like to thank Simone Gutt and Fred Van Oystaeyen for their kind hospitality.

\end{document}